\documentclass[11pt,a4paper]{article}
\usepackage[T1]{fontenc}
\usepackage{lmodern,microtype,amsmath,amssymb,amsthm,booktabs,array,longtable}
\usepackage[margin=25mm]{geometry}
\usepackage{graphicx,xcolor,natbib,placeins}
\usepackage{algorithm}
\usepackage{algpseudocode}

\usepackage[colorlinks=true,linkcolor=red,citecolor=blue,urlcolor=blue!45!black]{hyperref}
\usepackage{needspace}
\newtheorem{definition}{Definition}
\newtheorem{theorem}{Theorem}
\newtheorem{proposition}[theorem]{Proposition}

\newtheorem{remark}{Remark}
\newcommand{\dist}{\operatorname{dist}}

\graphicspath{{figures/}}
\title{Compromise Optimization over Incompatible Efficient Sets Using a Minimax Distance Approach with Polyhedral Characterization and an Application to Fisheries Management}
\author{
 Mustapha Kaci\\
 Laboratory of Mathematics and Statistics and Their Applications, USTOMB.\\
		Djillali Liabes University of Sidi Bel Abbes, Sidi Bel Abbes, Algeria.\\
  {\color{blue}\texttt{kaci.mustapha.95@gmail.com} } and {\color{blue}\texttt{mustapha.kaci@univ-sba.dz}}} 
\date{}

\newtheorem{answer}{Answer}
\newtheorem{question}{Question}
 \usepackage{tikz}
\usetikzlibrary{arrows.meta,positioning}

\begin{document}
\maketitle

\begin{abstract}
The efficient sets associated with different decision-makers may have an empty intersection, even when their objectives are linear and the common feasible region is a compact polytope. In this paper, we consider a compromise criterion that minimizes the largest weighted distance from a feasible decision to the respective efficient sets of the decision-makers. First, we establish the existence of an optimal compromise and show that the formulation recovers the common efficient set whenever it is nonempty. We then derive a lower bound based on the separation of efficient sets and obtain an exact formulation in terms of finitely many linear programs by exploiting their polyhedral structure. Next, we examine the effect of approximating the efficient sets and establish a stability bound expressed in terms of the Hausdorff distance. Finally, we illustrate the applicability of the proposed approach through a numerical study of harvest allocation between two interconnected fishing zones.
\end{abstract}


\section{Introduction}
\label{sec:introduction}

In many decision-making problems, several decision-makers are required to agree on a common feasible solution, although their objectives may differ or even conflict. Each decision-maker has a set of efficient solutions determined by their own objectives. These solutions are characterized by the fact that improving one objective necessarily leads to the deterioration of at least one other. Even when the individual efficient sets are nonempty, their intersection may be empty, so that no feasible solution is simultaneously efficient for all decision-makers. This raises the question of how to determine a compromise that takes into account the different preferences while preserving feasibility. In this paper, we approach this question by considering the distances between a feasible decision and the efficient sets of the individual decision-makers.

The characterization of efficient solutions and the search for suitable compromises are two closely related problems in multiple objective optimization. The development of these problems and the main approaches proposed in the literature are reviewed by Ehrgott et al.~\cite{Ehrgott2026}. In multiobjective linear programming, the polyhedral structure of the feasible region plays an important role in the characterization of efficient solutions. In particular, the efficient set can be described in terms of faces of the feasible polyhedron. Among the early contributions, Yu and Zeleny~\cite{YuZeleny1975} investigated the geometric properties of nondominated solutions and developed a multicriteria simplex method. Their results showed that the study of efficient faces provides a description of the efficient set that is not restricted to efficient extreme points.

Several methods have subsequently been developed to exploit the facial structure of efficient sets in multiobjective linear programming. Armand~\cite{Armand1993} studied the identification of maximal efficient faces, while Say{\i}n~\cite{Sayin1996} proposed a top-down algorithm based on facial decomposition and a linear programming test for efficiency. This method examines faces of higher dimension before considering their subfaces, without requiring the preliminary enumeration of all efficient extreme points. Benson and Say{\i}n~\cite{BensonSayin1997} studied the representation of efficient sets and its use in decision-making. Later, Say{\i}n~\cite{Sayin2000} applied a top-down search procedure to the optimization of a linear function over an efficient set. Tohidi and Hassasi~\cite{TohidiHassasi2018} extended this line of research by incorporating adjacency relations into a local search for maximal efficient faces. These methods take advantage of the polyhedral structure of the feasible region to identify and represent efficient sets as finite unions of faces, even when these sets are not convex.   In a different direction, Benson~\cite{Benson1998} developed an outer approximation algorithm for generating efficient extreme points of the attainable objective set, working in objective space rather than directly in the space of feasible decisions.

The search for compromises between different systems of objectives has also been investigated in multilevel programming, where the decision-makers are organized according to a hierarchy of priorities. Kaci and Radjef~\cite{KaciRadjef2022} studied the set of feasible solutions that are simultaneously nondominated with respect to the objective systems of all levels. Their approach identifies the set of possible compromises before taking the hierarchical structure into account in the selection process. Subsequently, Kaci and Radjef~\cite{KaciRadjef2023} developed an adaptive method for determining a compromise solution. These works provide a basis for studying common efficiency in the presence of several decision-makers. However, the existence of a common efficient solution is not always guaranteed, since the individual efficient sets may have an empty intersection. In this paper, we consider this situation by replacing the requirement of simultaneous efficiency with a distance-based criterion. More precisely, we seek a feasible solution that minimizes the largest weighted distance to the individual efficient sets.

The proposed formulation is related to the smallest intersecting ball problem, studied by Nam, Nguyen and Salinas~\cite{NamNguyenSalinas2012} using convex analysis. Related existence and uniqueness results were established by Mordukhovich, Nam and Villalobos~\cite{MordukhovichNamVillalobos2013}, while Nam, Villalobos and An~\cite{NamVillalobosAn2012} considered extensions involving minimal-time functions. In our problem, the compromise is restricted to the feasible polytope, and the reference sets correspond to the individual efficient sets of the decision-makers. For allocation problems with local capacities and an aggregate resource constraint, we obtain an explicit characterization of optimal compromises based on a scalar threshold, related to the resource-allocation methods surveyed by Patriksson~\cite{Patriksson2008}. As an application, we consider harvest allocation between two interconnected fishing zones. Elbetch, Moussaoui and Auger~\cite{Elbetch2025} studied maximum economic yield in heterogeneous fishing patches with fast migration. Here, we adopt a fixed-period linear biomass model, using publicly available aggregate data to set the numerical scale and specified assumptions to describe the spatial allocation and local management constraints.

The remainder of the paper is organized as follows. Section~\ref{sec:mlmolp-formulation} introduces the mathematical formulation of the problem and the feasible region. In Section~\ref{sec:preliminaries}, we recall the definitions of individual efficient sets and the set of possible compromises. Section~\ref{sec:mainresults} presents the proposed distance-based criterion and establishes its main properties, including existence of optimal solutions, recovery of the common efficient set and separation bounds. Section~\ref{sec:lp-formulation} deals with the finite linear programming formulation and the progressive generation of efficient faces, while Section~\ref{secapp} examines the effect of approximating the reference sets. Section~\ref{sec:application} presents the harvest allocation example. The results are discussed in Section~\ref{sec:discussion}, and the paper concludes in Section~\ref{sec:conclusion}.

\section{Problem formulation}
\label{sec:mlmolp-formulation}

We consider a multi-level multi-objective linear programming problem involving $P\geq 2$ decision-makers, denoted by $\mathrm{DM}_p$, $p=1,\ldots,P$, and organized according to a hierarchical structure. The first decision-maker, $\mathrm{DM}_1$, occupies the highest level, while $\mathrm{DM}_P$ is at the lowest level. Each decision-maker $\mathrm{DM}_p$ controls a vector of variables $x_p\in\mathbb{R}^{n_p}$. The complete decision vector is given by
\[
x=(x_1^\top,\ldots,x_P^\top)^\top\in\mathbb{R}^n,
\qquad n=\sum_{p=1}^{P}n_p.
\]

At each level $p$, the decision-maker seeks to maximize $k_p\geq 1$ linear objective functions, represented by
\[
F_p(x)=C_px=
\bigl(c_{p1}^\top x,\ldots,c_{pk_p}^\top x\bigr)^\top
\in\mathbb{R}^{k_p},
\qquad C_p\in\mathbb{R}^{k_p\times n}.
\]
The objectives are considered in the componentwise maximization sense. All decision-makers are subject to the same feasible region
\begin{equation}\label{eq:common-feasible-region}
S=\{x\in\mathbb{R}^n:Ax\leq b,\ x\geq 0\},
\end{equation}
where $A\in\mathbb{R}^{m\times n}$ and $b\in\mathbb{R}^m$. We assume throughout that $S$ is nonempty and compact. Consequently, $S$ is a bounded convex polyhedron.

The hierarchical problem can be expressed schematically as
\begin{equation}\label{eq:mlmolp}
\begin{aligned}
&\text{Level }1:\quad \max_{x_1}\ F_1(x),
&&\text{subject to the decisions at lower levels},\\
&\text{Level }2:\quad \max_{x_2}\ F_2(x),
&&\text{subject to the decisions at lower levels},\\
&\hspace{3.6em}\vdots\\
&\text{Level }P:\quad \max_{x_P}\ F_P(x),
&&\text{subject to }x\in S.
\end{aligned}
\end{equation}

\subsection{Hierarchical resolution procedure}
\label{sec:hierarchical-resolution}

The resolution of a multilevel problem follows the hierarchy of decision-makers. As described in~\cite{Kaci2025}, each level solves its optimization problem and restricts the feasible region for the following level. The procedure consists of three main steps:

\begin{enumerate}
    \item \textbf{First level:} The upper-level decision-maker solves
    \[
    \max_{x\in S} f_1(x).
    \]
    Its solution is used to define a reduced feasible region $S_1\subseteq S$.

    \item \textbf{Intermediate levels:} Each decision-maker $\mathrm{DM}_p$, $p=2,\ldots,P-1$, solves
    \[
    \max_{x\in S_{p-1}} f_p(x),
    \]
    and introduces additional restrictions, giving
    $S_p\subseteq S_{p-1}$.

    \item \textbf{Final level:} The last decision-maker solves
    \[
    \max_{x\in S_{P-1}} f_P(x).
    \]
    An optimal solution of this problem is retained as the final compromise.
\end{enumerate}

Thus, the procedure generates a sequence of nested feasible regions,
\[
S\supseteq S_1\supseteq\cdots\supseteq S_{P-1},
\]
which are assumed to remain nonempty throughout the process.

\begin{figure}[H]
\centering
\begin{tikzpicture}[
    >=Latex,
    node distance=7mm,
    every node/.style={font=\small},
    box/.style={
        draw=black!70,
        rounded corners=6pt,
        very thick,
        align=center,
        minimum height=1cm,
        text width=8.7cm,
        inner sep=6pt
    },
    start/.style={
        box,
        fill=gray!10
    },
    level/.style={
        box,
        fill=blue!8
    },
    restrict/.style={
        box,
        fill=teal!8
    },
    final/.style={
        box,
        fill=orange!10
    },
    arrow/.style={
        ->,
        very thick,
        draw=black!70
    }
]

\node[start] (s0) {%
\textbf{Initial feasible region}\\
$S$
};

\node[level, below=of s0] (l1) {%
\textbf{Level 1}\\
Solve $\displaystyle \max_{x\in S} f_1(x)$
};

\node[restrict, below=of l1] (r1) {%
\textbf{Restriction}\\
Obtain a reduced feasible region $S_1\subseteq S$
};

\node[level, below=of r1] (lp) {%
\textbf{Intermediate levels \;$p=2,\dots,P-1$}\\
Solve $\displaystyle \max_{x\in S_{p-1}} f_p(x)$
};

\node[restrict, below=of lp] (rp) {%
\textbf{Further restriction}\\
Obtain $S_p\subseteq S_{p-1}$
};

\node[below=2mm of rp] (dots) {$\vdots$};

\node[level, below=2mm of dots] (lP) {%
\textbf{Final level \;$P$}\\
Solve $\displaystyle \max_{x\in S_{P-1}} f_P(x)$
};

\node[final, below=of lP] (xc) {%
\textbf{Final compromise}\\
An optimal solution is retained as the final compromise
};

\draw[arrow] (s0) -- (l1);
\draw[arrow] (l1) -- (r1);
\draw[arrow] (r1) -- (lp);
\draw[arrow] (lp) -- (rp);
\draw[arrow] (rp) -- (dots);
\draw[arrow] (dots) -- (lP);
\draw[arrow] (lP) -- (xc);

\end{tikzpicture}
\caption{Hierarchical resolution procedure through successive reductions of the feasible region.}
\label{fig:hierarchical-resolution}
\end{figure}
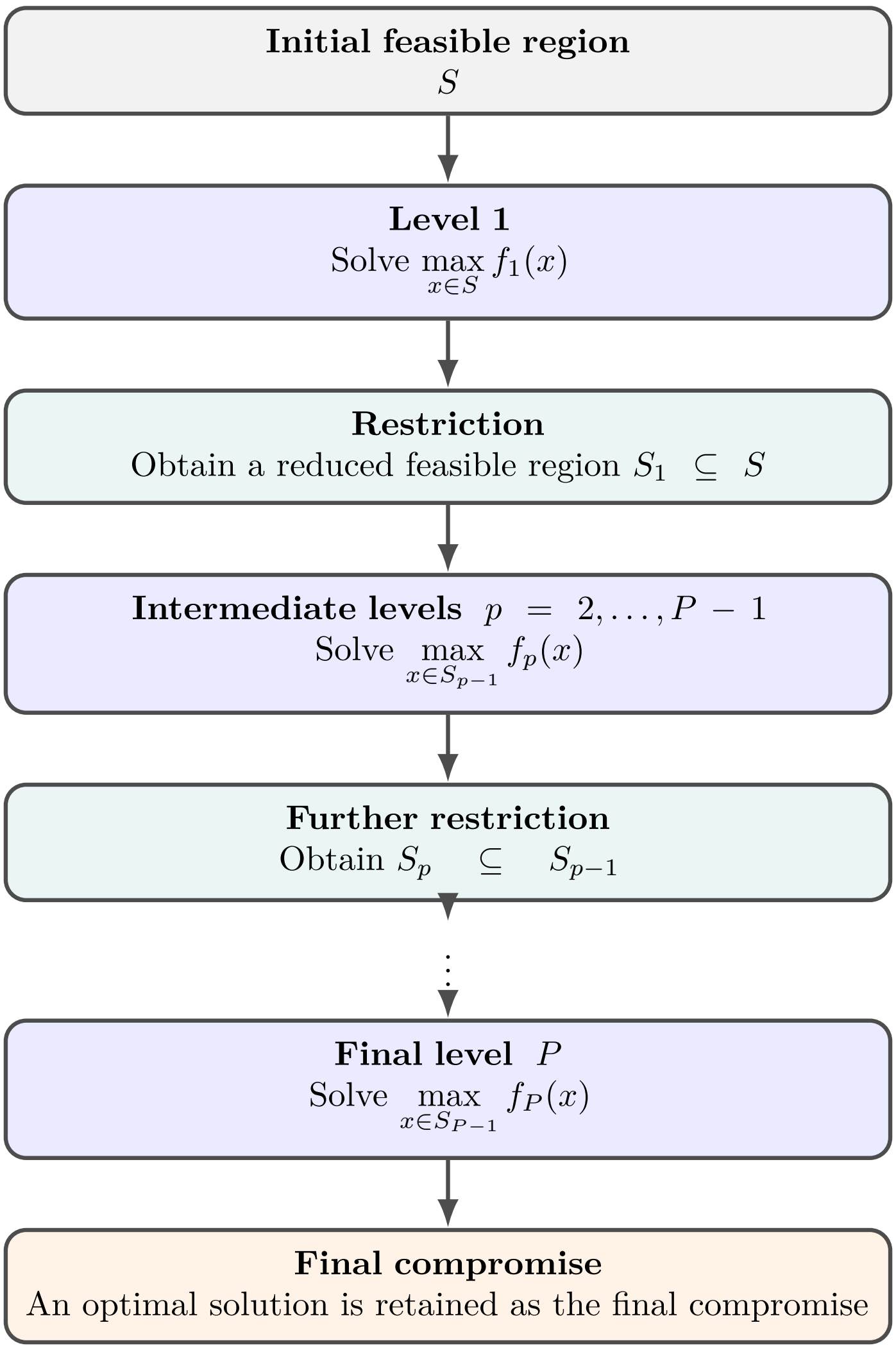
\section{Preliminaries}
\label{sec:preliminaries}

For each decision-maker $\mathrm{DM}_p$, $p=1,\ldots,P$, we consider the multi-objective linear programming problem
\begin{equation}\label{eq:level-molp}
    \max\{F_p(x):x\in S\},
\end{equation}
where the objectives are maximized with respect to the componentwise order in $\mathbb{R}^{k_p}$. The corresponding attainable objective set is
\[
Z_p=\{F_p(x):x\in S\}.
\]

\begin{definition}[Dominance and non-dominance]Let $x,y\in S$, then
\begin{enumerate}
    \item  .  We say that $y$ \emph{dominates} $x$ at level $p$ if
\[
F_p(y)\geq F_p(x),
\qquad F_p(y)\neq F_p(x),
\]
where the inequality is understood componentwise.
\item
The point $x$ is said to be \emph{non-dominated} at level $p$ if there is no $y\in S$ that dominates it. The set of all such points is denoted by
\begin{equation}\label{eq:level-nondominated}
E_p=
\left\{
x\in S:
\nexists\,y\in S,\;
F_p(y)\geq F_p(x),\;
F_p(y)\neq F_p(x)
\right\}\neq\emptyset.
\end{equation}
\end{enumerate}
\end{definition}

\begin{definition}[Set of all possible compromises, \cite{KaciRadjef2022}]
The \emph{set of all possible compromises} of the multilevel multi-objective linear programming problem is defined by
\begin{equation}\label{eq:all-possible-compromises}
E=\bigcap_{p=1}^{P}E_p.
\end{equation}
Thus, $E$ consists of all feasible solutions that are non-dominated with respect to the objective system of every decision-maker.
\end{definition}

Although each individual efficient set $E_p$ is nonempty, their intersection $E$ may be empty. This raises the following question:

\begin{question}
How can a compromise solution be determined for the multilevel problem when $E=\varnothing$?
\end{question}

\begin{answer}
We propose a distance-based compromise criterion that seeks a feasible solution as close as possible to the individual efficient sets, without requiring their intersection to be nonempty.
\end{answer}

For convenience, the main notation used throughout the paper is summarized in Table~\ref{tab:notation}.

\begin{longtable}{p{0.19\textwidth}p{0.73\textwidth}}
\caption{Main notation used throughout the paper.}
\label{tab:notation}\\
\toprule
Symbol & Meaning\\
\midrule
\endfirsthead

\toprule
Symbol & Meaning\\
\midrule
\endhead

$P$, $p$
& Number of decision-makers and decision-maker index.\\

$n$, $x$
& Total number of decision variables and complete decision vector.\\

$n_p$, $x_p$
& Number of variables and decision vector controlled by $\mathrm{DM}_p$.\\

$S$
& Nonempty compact polyhedral feasible region.\\

$k_p$, $C_p$
& Number of objectives and objective coefficient matrix of $\mathrm{DM}_p$.\\

$F_p(x)$
& Vector of linear objectives associated with $\mathrm{DM}_p$.\\

$E_p$
& Set of non-dominated feasible solutions for $\mathrm{DM}_p$.\\

$E$
& Intersection of the individual efficient sets.\\

$\sigma_i$
& Positive scaling factor associated with decision variable $i$.\\

$\|\cdot\|_\sigma$
& Scaled maximum norm in decision space.\\

$\alpha_p$
& Positive weight assigned to decision-maker $\mathrm{DM}_p$.\\

$d_p(x)$
& Distance from a feasible point $x$ to the efficient set $E_p$.\\

$\Phi(x)$
& Maximum weighted distance from $x$ to the individual efficient sets.\\

$\tau^*$
& Minimum value of the compromise criterion over $S$.\\

$X^*$
& Set of optimal compromise solutions.\\

$\delta_{pq}$
& Minimum distance between the efficient sets $E_p$ and $E_q$.\\

$Q_{pr}$
& Efficient face indexed by $r$ for decision-maker $\mathrm{DM}_p$.\\

$R_p$
& Number of retained efficient faces associated with $\mathrm{DM}_p$.\\

\bottomrule
\end{longtable}

%%%%%%%%%%%%%%%%%%%%%%%%%%%%%%%%%%%%%%%%%%%%%%%%%%%%%%%%%%%%%%%%%%%%%%%%%%%%%%%%%%%%%%%%%%%%%%%%%%%%%%%%%%%%%%%%%%%%%%%%%%%%%%%%%%%%%%%%%%%%%%%%%%%%%%%%%%%%%%%%%%%%%%%%%%%%%%%%%%%%%%%%%%%%%%%%%%%%%%%%%%%%%%%%%%%%%%%%%%%%%%%%%%%%%%%%%%%%%%%%%%%%%%%%%%%%%%%%%%%%%%%%%%%%%%%%%%%%%%%%%%%%%%%%%%%%%%%%%%%%%%%%%%%%%%%%%%%%%%%%%%%%%%%%%%%%%%%%%%%
\section{Main results}\label{sec:mainresults}
\subsection{Distance to the individual efficient sets}

When the individual efficient sets $E_p$ have no common point, we seek a feasible decision whose distance to each of these sets is controlled. The distance is measured in decision space and represents the changes in the decision variables required to reach an efficient solution, rather than differences in objective values.

To account for differences in units and magnitudes, let $\sigma_i>0$, $i=1,\ldots,n$, and consider the scaled maximum norm
\[
\|z\|_\sigma
=\max_{1\leq i\leq n}\frac{|z_i|}{\sigma_i}.
\]
\begin{remark}$\;$\\

\begin{itemize}
    \item For each nonconstant coordinate, $\sigma_i=\max_{x\in S}x_i-\min_{x\in S}x_i.$
    \item Alternatively, $\sigma_i$ may represent a prescribed variation in the corresponding decision variable. Coordinates that are constant on $S$ may be omitted or assigned any positive scale.
\end{itemize}
\end{remark}

Since $S$ is a nonempty compact polytope and the objectives are linear, each efficient set $E_p$ is nonempty and can be represented as a finite union of compact faces~\cite{YuZeleny1975}. Thus, the distance from a feasible point $x$ to $E_p$ is attained and is defined by
\[
d_p(x)=\operatorname{dist}_\sigma(x,E_p)
=\min_{y\in E_p}\|x-y\|_\sigma.
\]
In particular,
\[
d_p(x)=0 \quad\Longleftrightarrow\quad x\in E_p.
\]

To account for the relative importance of the decision-makers, let $\alpha_p>0$, $p=1,\ldots,P$. We define the compromise problem as
\begin{equation}
\label{model}
\tau^*=\min_{x\in S}\Phi(x),
\qquad
\Phi(x)=\max_{1\leq p\leq P}\alpha_p d_p(x).
\end{equation}

The criterion $\Phi$ measures the largest weighted distance to the individual efficient sets. Hence, $\tau^*$ represents the smallest common tolerance that can be achieved by a feasible decision. The weights $\alpha_p$ determine the relative importance of the decision-makers without imposing a hierarchical order.

\begin{remark}
The distances $d_p(x)$ are measured in decision space, not in objective space. They quantify the changes in decision variables needed to reach an efficient solution, rather than the corresponding loss in objective values. 
\end{remark}

For $\tau\geq 0$, let $K(\tau)$ denote the set of feasible decisions whose weighted distances to all individual efficient sets do not exceed $\tau$:
\[
K(\tau)=\left\{x\in S:
\alpha_p d_p(x)\leq\tau,\quad p=1,\ldots,P
\right\}.
\]
Equivalently,
\[
K(\tau)=S\cap\bigcap_{p=1}^{P}
\left(E_p+\frac{\tau}{\alpha_p}\mathbb{B}_\sigma\right),
\]
where $\mathbb{B}_\sigma$ is the unit ball of the scaled maximum norm. In particular, at zero tolerance, we recover the set of all possible compromises:
\[
K(0)=\bigcap_{p=1}^{P}E_p=E.
\]

\begin{definition}[$\tau$-efficient compromise solution]
Let $\tau\geq 0$. A feasible point $x\in S$ is called a
$\tau$-efficient compromise solution if
\[
\alpha_p d_p(x)\leq\tau,
\qquad p=1,\ldots,P,
\]
or, equivalently, if $x\in K(\tau)$.

The smallest admissible tolerance is defined by
\[
\tau^*=\min\{\tau\geq 0:K(\tau)\neq\varnothing\}
      =\min_{x\in S}\Phi(x).
\]
A $\tau^*$-efficient compromise solution is called an
\emph{optimal compromise solution}. The set of all such
solutions is
\begin{equation}
\label{eq:optimal-compromise-set}
X^*=\operatorname*{arg\,min}_{x\in S}\Phi(x)
   =K(\tau^*).
\end{equation}

In particular,
\begin{itemize}
    \item If $E\neq\varnothing$, then $\tau^*=0$ and $X^*=E$.
    \item If $E=\varnothing$, then $\tau^*>0$, and $X^*$
    consists of the feasible decisions attaining the
    smallest possible maximum weighted distance to
    the individual efficient sets.
\end{itemize}
\end{definition}

%%%%%%%%%%%%%%%%%%%%%%%%%%%%%%%%%%%%%%%%%%%%%%%%%%%%%%%%%%%%%%%%%%%%%%%%%%%%%%%%%%%%%%%%%%%%%%%%%%%%%%%%%%%%%%%%%%%%%%%%%%%%%%%%%%%%%%%%%%%%%%%%%%%%%%%%%%%%%%%%%%%%%%%%%%%%%%%%%%%%%%%%%%%%%%%%%%%%%%%%%%%%%%%%%%%%%%%%%%%%%%%%%%%%%%%%%%%%%%%%%%%%%%%%%%%%%%%%%%%%%%%%%%%%%%%%%%%%%%%%%%%%%%%%%%%%%%%%%%%%%%%%%%%%%%%%%%%%%%%%%%%%%%%%%%%%%%%%%%%
\subsection{Geometric Characterization of Optimal Compromises}

\begin{theorem}[Existence and exact recovery of common efficiency]
\label{thm:existence}
Let $S\subset\mathbb{R}^n$ be a nonempty compact
polytope, and let $E_p\subseteq S$ be nonempty
closed sets for $p=1,\ldots,P$. Assume that
$\sigma_i>0$ for $i=1,\ldots,n$ and
$\alpha_p>0$ for $p=1,\ldots,P$.
Then the following assertions hold:

\begin{enumerate}
    \item \textbf{Existence.}
    Problem~\eqref{model} admits an optimal solution, and
    $X^*=\operatorname*{arg\,min}_{x\in S}\Phi(x)$
    is nonempty and compact.

    \item \textbf{Exact recovery.}
    \begin{equation}
    \label{eq:exact-recovery}
    \tau^*=0
    \iff E=\bigcap_{p=1}^{P}E_p\neq\varnothing.
    \end{equation}
    Moreover, if $E\neq\varnothing$, then
    \begin{equation}
    \label{eq:exact-minimizers}
    X^*=E.
    \end{equation}

    \item \textbf{Tolerance characterization.}
    For every $\tau\geq 0$,
    \begin{equation}
    \label{eq:radius-threshold}
    K(\tau)\neq\varnothing
    \iff \tau\geq\tau^*,
    \end{equation}
    and
    \begin{equation}
    \label{eq:optimal-level-set}
    X^*=K(\tau^*).
    \end{equation}
\end{enumerate}

\end{theorem}

\begin{proof}$\,$\\
\begin{enumerate}
    \item 

We first show that the compromise objective
$\Phi$ is continuous on $S$.

Fix $p\in\{1,\ldots,P\}$. Since $E_p$ is a
nonempty closed subset of the compact set $S$,
it is compact. Consequently, for every $x\in S$,
there exists $y_p(x)\in E_p$ such that
\[
d_p(x)=\|x-y_p(x)\|_\sigma.
\]

Let $x,z\in S$. By the triangle inequality,
\[
\|x-y\|_\sigma
\leq
\|x-z\|_\sigma+\|z-y\|_\sigma
\qquad \text{for every }y\in E_p.
\]
Taking the minimum over $y\in E_p$ gives
\begin{equation}\label{a}
    d_p(x)\leq\|x-z\|_\sigma+d_p(z).
\end{equation}

Interchanging $x$ and $z$, we similarly obtain
\begin{equation}\label{aa}
d_p(z)\leq\|x-z\|_\sigma+d_p(x).
\end{equation}
Combining these inequalities (\ref{a})-(\ref{aa}), we get
\begin{equation}
\label{eq:distance-lipschitz}
|d_p(x)-d_p(z)|
\leq\|x-z\|_\sigma.
\end{equation}
Thus $d_p$ is Lipschitz continuous on $S$.

Since $\alpha_p>0$, each function
$x\mapsto\alpha_p d_p(x)$ is continuous.
The compromise objective
\[
\Phi(x)=\max_{1\leq p\leq P}\alpha_p d_p(x)
\]
is therefore continuous as the pointwise
maximum of finitely many continuous functions.

As $S$ is nonempty and compact, the
Weierstrass theorem guarantees the existence
of $x^*\in S$ such that
\[
\Phi(x^*)=\min_{x\in S}\Phi(x)=\tau^*.
\]
Hence $X^*$ is nonempty. Moreover,
\[
X^*=\{x\in S:\Phi(x)=\tau^*\}
\]
is closed in $S$ by continuity of $\Phi$.
Since $S$ is compact, $X^*$ is compact.

\item
We first establish the equivalence
\eqref{eq:exact-recovery}.

Suppose that $E\neq\varnothing$ and choose
$x\in E$. By the definition of the common
efficient set,
\[
x\in E_p
\qquad\text{for every }p=1,\ldots,P.
\]
Consequently, $d_p(x)=0$ for every $p$, and
hence
\[
\Phi(x)=\max_{1\leq p\leq P}\alpha_p d_p(x)=0.
\]
Since $\Phi$ is nonnegative on $S$, it follows
that $\tau^*=0$.

Conversely, suppose that $\tau^*=0$.
By assertion (i), there exists $x^*\in X^*$
such that
\[
0=\Phi(x^*)
=\max_{1\leq p\leq P}\alpha_p d_p(x^*).
\]
All terms in this maximum are nonnegative.
Since the weights $\alpha_p$ are strictly
positive, we necessarily have
\[
d_p(x^*)=0
\qquad\text{for every }p=1,\ldots,P.
\]
Because $E_p$ is closed and nonempty,
the distance to $E_p$ vanishes precisely
at points belonging to $E_p$. Thus
\[
x^*\in E_p
\qquad\text{for every }p=1,\ldots,P,
\]
which implies $x^*\in E$. In particular,
$E\neq\varnothing$, proving
\eqref{eq:exact-recovery}.

We now prove \eqref{eq:exact-minimizers}.
Assume that $E\neq\varnothing$. By the
preceding equivalence, $\tau^*=0$.

For every $x\in E$, all individual distances
vanish, so $\Phi(x)=0=\tau^*$. Therefore,
\[
E\subseteq X^*.
\]
Conversely, if $x\in X^*$, then
$\Phi(x)=\tau^*=0$. The same argument as
above shows that $d_p(x)=0$ for every $p$,
and hence $x\in E$. Therefore,
\[
X^*\subseteq E.
\]
The two inclusions establish $X^*=E$.

\item
Let $\tau\geq0$. By the definition of
$K(\tau)$ and the positivity of the weights,
\begin{align*}
x\in K(\tau)
&\iff
d_p(x)\leq\frac{\tau}{\alpha_p}
\quad\text{for every }p\\
&\iff
\alpha_p d_p(x)\leq\tau
\quad\text{for every }p\\
&\iff
\max_{1\leq p\leq P}\alpha_p d_p(x)
\leq\tau\\
&\iff
\Phi(x)\leq\tau.
\end{align*}
Consequently,
\begin{equation}
\label{eq:K-sublevel}
K(\tau)=\{x\in S:\Phi(x)\leq\tau\}.
\end{equation}

Suppose first that $K(\tau)\neq\varnothing$.
Then there exists $x\in S$ such that
$\Phi(x)\leq\tau$. Since
\[
\tau^*=\min_{z\in S}\Phi(z)\leq\Phi(x),
\]
we obtain $\tau^*\leq\tau$.

Conversely, suppose that $\tau\geq\tau^*$.
By assertion (i), there exists an optimal
solution $x^*\in X^*$ satisfying
\[
\Phi(x^*)=\tau^*\leq\tau.
\]
Equation \eqref{eq:K-sublevel} then implies
$x^*\in K(\tau)$, so $K(\tau)$ is nonempty.
This proves \eqref{eq:radius-threshold}.

Finally, setting $\tau=\tau^*$ in
\eqref{eq:K-sublevel} gives
\[
K(\tau^*)
=\{x\in S:\Phi(x)\leq\tau^*\}.
\]
By the definition of $\tau^*$, every
$x\in S$ satisfies $\Phi(x)\geq\tau^*$.
It follows that
\[
K(\tau^*)
=\{x\in S:\Phi(x)=\tau^*\}
=X^*,
\]
which establishes \eqref{eq:optimal-level-set}
and completes the proof.
\end{enumerate}
\end{proof}

\begin{proposition}[Lower Bound on the Optimal Compromise Radius]\label{n}
For $p\ne q$, define
$\delta_{pq}=\min\{\|u-v\|_\sigma:u\in E_p,v\in E_q\}$.
Then
\begin{equation}
\tau^*\ge\max_{p<q}\frac{\alpha_p\alpha_q}{\alpha_p+\alpha_q}\delta_{pq}.
\label{bound}
\end{equation}
\end{proposition}

\begin{proof}
Fix $1\leq p<q\leq P$. Since $E_p$ and $E_q$ are nonempty and compact, the distance $\delta_{pq}$ is attained.

For any $x\in S$, choose $u\in E_p$ and $v\in E_q$ such that
\[
d_p(x)=\|x-u\|_\sigma,
\qquad
d_q(x)=\|x-v\|_\sigma.
\]
By the triangle inequality,
\begin{equation}
\label{eq:pairwise-separation}
\delta_{pq}
\leq \|u-v\|_\sigma
\leq d_p(x)+d_q(x).
\end{equation}
On the other hand, the definition of $\Phi$ gives
\[
d_p(x)\leq\frac{\Phi(x)}{\alpha_p},
\qquad
d_q(x)\leq\frac{\Phi(x)}{\alpha_q}.
\]
Consequently,
\[
\delta_{pq}
\leq
\left(\frac{1}{\alpha_p}
+\frac{1}{\alpha_q}\right)\Phi(x),
\]
which implies
\[
\Phi(x)\geq
\frac{\alpha_p\alpha_q}{\alpha_p+\alpha_q}
\,\delta_{pq}.
\]
Taking the minimum over $x\in S$, we obtain
\[
\tau^*\geq
\frac{\alpha_p\alpha_q}{\alpha_p+\alpha_q}
\,\delta_{pq}.
\]
Since the inequality holds for every pair $p<q$, taking the maximum over all such pairs proves~\eqref{bound}.
\end{proof}

The bound \eqref{bound} provides a lower estimate of the
optimal compromise radius in terms of the pairwise distances
between the individual efficient sets. In particular, if
$E_p\cap E_q=\varnothing$ for some $p\neq q$, then
$\delta_{pq}>0$, and hence $\tau^*>0$.

Any feasible point whose objective value attains the bound
\eqref{bound} is necessarily optimal. However, the bound
is not generally sharp when $P\geq 3$. Indeed, all pairwise
distances may vanish even though the common intersection
$\bigcap_{p=1}^{P}E_p$ is empty.

\begin{theorem}[An exact construction for two decision makers]\label{two}
Suppose $P=2$. Let $(u^*,v^*)\in E_1\times E_2$ attain $\delta_{12}$. Then
\begin{equation}\label{eq:two-solution}
\tau^*=\frac{\alpha_1\alpha_2}{\alpha_1+\alpha_2}\delta_{12},\qquad
\bar x=\frac{\alpha_1u^*+\alpha_2v^*}{\alpha_1+\alpha_2}\in X^*.
\end{equation}
The efficient sets themselves need not be convex.
\end{theorem}

\begin{proof}
Since $E_1$ and $E_2$ are nonempty and compact,
a pair $(u^*,v^*)$ attaining $\delta_{12}$ exists.
The positivity of the weights and the convexity
of $S$ ensure that the point $\bar x$ defined
in \eqref{eq:two-solution} belongs to $S$.

We first establish an upper bound for $\tau^*$.
Using $u^*\in E_1$ and the definition of $d_1$,
we obtain
\begin{align*}
\alpha_1d_1(\bar x)
&\leq \alpha_1\|\bar x-u^*\|_\sigma\\
&=\frac{\alpha_1\alpha_2}{\alpha_1+\alpha_2}
  \|v^*-u^*\|_\sigma\\
&=\frac{\alpha_1\alpha_2}{\alpha_1+\alpha_2}
  \delta_{12}.
\end{align*}
Similarly, since $v^*\in E_2$,
\[
\alpha_2d_2(\bar x)
\leq
\frac{\alpha_1\alpha_2}{\alpha_1+\alpha_2}
\,\delta_{12}.
\]
Consequently,
\begin{equation}
\label{eq:two-upper}
\tau^*
\leq \Phi(\bar x)
\leq
\frac{\alpha_1\alpha_2}{\alpha_1+\alpha_2}
\,\delta_{12}.
\end{equation}

On the other hand, Proposition~\ref{n}
gives
\begin{equation}
\label{eq:two-lower}
\tau^*
\geq
\frac{\alpha_1\alpha_2}{\alpha_1+\alpha_2}
\,\delta_{12}.
\end{equation}
Combining \eqref{eq:two-upper} and
\eqref{eq:two-lower} establishes
\eqref{eq:two-solution}. Moreover, equality
throughout \eqref{eq:two-upper} yields
$\Phi(\bar x)=\tau^*$, proving that
$\bar x\in X^*$.

The argument relies only on the existence
of a nearest pair in $E_1\times E_2$ and
the convexity of the feasible set $S$.
Thus no convexity assumption on $E_1$
or $E_2$ is needed.
\end{proof}

Theorem~\ref{two} provides an explicit construction of an
optimal compromise for two decision makers from a pair of
nearest points in their respective efficient sets. This
construction does not, however, characterize the entire
set of minimizers, which need not be a singleton, particularly
under the maximum norm. For $P\geq 3$, the pairwise lower
bound \eqref{bound} remains valid but is not necessarily attained.

\subsection{An Exact Linear Programming Formulation}
\label{sec:lp-formulation}

The compromise problem can be reformulated as a finite
family of linear programs when the individual efficient
sets admit explicit polyhedral representations. The
reformulation relies on the fact that each efficient set
is a finite union of faces of the feasible polytope.
Although these sets need not be convex, each individual
face is a convex polytope and can be described by linear
constraints.

For each decision maker $p\in\{1,\ldots,P\}$, suppose
that the efficient set $E_p$ is represented as
\begin{equation}
\label{eq:efficient-face-decomposition}
E_p=\bigcup_{r=1}^{R_p}Q_{pr},
\end{equation}
where $R_p$ denotes the number of faces retained in
the representation and $Q_{pr}$ is the $r$-th efficient
face associated with decision maker $p$.

Each face is assumed to be nonempty and is described by
\begin{equation}
\label{eq:face-representation}
Q_{pr}
=
\left\{
y\in S:
H_{pr}y=h_{pr}
\right\},
\end{equation}
where
$H_{pr}\in\mathbb{R}^{m_{pr}\times n}$ is a matrix
whose rows specify the linear equalities defining
the face, and
$h_{pr}\in\mathbb{R}^{m_{pr}}$ is the corresponding
right-hand-side vector. Here, $m_{pr}$ denotes the
number of equalities used to describe $Q_{pr}$.
The inequalities defining $S$ remain imposed through
the condition $y\in S$.

When $S=\{x\in\mathbb{R}^n:Ax\leq b\}$, the equalities
in \eqref{eq:face-representation} may be obtained
from inequalities of $Ax\leq b$ that are active
throughout the face. In particular, if
$I_{pr}$ is an index set of such inequalities,
one may take
\[
H_{pr}=A_{I_{pr}},
\qquad
h_{pr}=b_{I_{pr}}.
\]
The representation need not be irredundant.
However, it is sufficient to retain the
inclusion-maximal efficient faces, since every
efficient face contained in another retained
face is redundant in \eqref{eq:efficient-face-decomposition}.

\subsection{Linear programming reformulation}

Define the finite index set
\[
\mathcal{R}
=
\prod_{p=1}^{P}\{1,\ldots,R_p\}.
\]
An element
$\boldsymbol r=(r_1,\ldots,r_P)\in\mathcal{R}$
specifies one efficient face $Q_{p r_p}$
for each decision maker $p$.

For a fixed tuple $\boldsymbol r\in\mathcal{R}$,
consider the following linear program:
\begin{equation}
\label{lp}
\begin{aligned}
\tau_{\boldsymbol r}
=\min_{x,y_1,\ldots,y_P,t}\quad &t\\
\text{subject to}\quad
&Ax\leq b,\\
&Ay_p\leq b,
&&p=1,\ldots,P,\\
&H_{p r_p}y_p=h_{p r_p},
&&p=1,\ldots,P,\\
&-\frac{\sigma_i}{\alpha_p}t
\leq x_i-y_{pi}
\leq\frac{\sigma_i}{\alpha_p}t,
&&\substack{i=1,\ldots,n,\\p=1,\ldots,P},\\
&t\geq0.
\end{aligned}
\end{equation}

Here, $x\in\mathbb{R}^n$ is the compromise decision,
$y_p\in\mathbb{R}^n$ is a reference decision
restricted to the selected face $Q_{p r_p}$,
and $t\in\mathbb{R}$ is an upper bound on the
weighted distances from $x$ to these reference
decisions.

Indeed, the coordinate inequalities in \eqref{lp}
are equivalent to
\[
\alpha_p\|x-y_p\|_\sigma\leq t,
\qquad p=1,\ldots,P.
\]
Thus, for a fixed selection of efficient faces,
the program minimizes the largest weighted
distance between the compromise decision and
one reference point chosen on each face.

\begin{theorem}[Exact finite linear programming reformulation]
\label{thm:finite-lp}
Suppose that $S$ is a nonempty compact polytope,
that $\alpha_p>0$ and $\sigma_i>0$, and that
the efficient sets admit the finite representations
\eqref{eq:efficient-face-decomposition}--%
\eqref{eq:face-representation}.
Then the following assertions hold:
\begin{enumerate}
    \item For every
    $\boldsymbol r\in\mathcal{R}$,
    the linear program \eqref{lp} is feasible
    and attains a finite optimal value
    $\tau_{\boldsymbol r}$.

    \item The optimal compromise radius satisfies
    \begin{equation}
    \label{eq:finite-lp-equivalence}
    \tau^*
    =
    \min_{\boldsymbol r\in\mathcal{R}}
    \tau_{\boldsymbol r}.
    \end{equation}

    \item If $\boldsymbol r^*\in\mathcal{R}$
    attains the minimum in
    \eqref{eq:finite-lp-equivalence}, then
    the $x$-component of any optimal solution
    of \eqref{lp} associated with
    $\boldsymbol r^*$ belongs to $X^*$.
\end{enumerate}
\end{theorem}

\begin{proof}$\;$\\
\begin{enumerate}
    \item 
Fix $\boldsymbol r\in\mathcal{R}$.
Since $S$ and the selected faces
$Q_{p r_p}$ are nonempty, choose
$x^0\in S$ and $y_p^0\in Q_{p r_p}$
for every $p$.

Set
\[
t^0=
\max_{\substack{1\leq p\leq P\\1\leq i\leq n}}
\frac{\alpha_p}{\sigma_i}
|x_i^0-y_{pi}^0|.
\]
Then $(x^0,y_1^0,\ldots,y_P^0,t^0)$
satisfies all constraints of \eqref{lp}.
Hence the program is feasible and
$\tau_{\boldsymbol r}\leq t^0$.

To establish attainment, it suffices to
restrict attention to feasible points with
$0\leq t\leq t^0$. Since $S$ and all
$Q_{p r_p}$ are compact, and the constraints
are closed, this restricted feasible set
is nonempty and compact. The continuous
objective $t$ therefore attains its minimum.

\medskip
\item
We establish both inequalities in
\eqref{eq:finite-lp-equivalence}.

First, consider any
$\boldsymbol r\in\mathcal{R}$ and let
$(x,y_1,\ldots,y_P,t)$ be feasible
for \eqref{lp}. Since
$y_p\in Q_{p r_p}\subseteq E_p$, we have
\[
d_p(x)
\leq\|x-y_p\|_\sigma
\leq\frac{t}{\alpha_p},
\qquad p=1,\ldots,P.
\]
Consequently,
\[
\Phi(x)
=
\max_{1\leq p\leq P}\alpha_p d_p(x)
\leq t.
\]
Since $x\in S$, it follows that
$\tau^*\leq\Phi(x)\leq t$.
Taking the minimum over the feasible
solutions of \eqref{lp} yields
\[
\tau^*\leq\tau_{\boldsymbol r}.
\]
As this holds for every tuple,
\begin{equation}
\label{eq:lp-first-inequality}
\tau^*
\leq
\min_{\boldsymbol r\in\mathcal{R}}
\tau_{\boldsymbol r}.
\end{equation}

Conversely, let $x^*\in X^*$ be an
optimal compromise. For each decision
maker $p$, compactness of $E_p$ ensures
the existence of $y_p^*\in E_p$ such that
\[
d_p(x^*)=\|x^*-y_p^*\|_\sigma.
\]
By \eqref{eq:efficient-face-decomposition},
there exists an index
$r_p^*\in\{1,\ldots,R_p\}$ for which
$y_p^*\in Q_{p r_p^*}$.

Consider the tuple
$\boldsymbol r^*=(r_1^*,\ldots,r_P^*)$.
Since
\[
\alpha_p\|x^*-y_p^*\|_\sigma
=
\alpha_p d_p(x^*)
\leq\Phi(x^*)=\tau^*,
\]
the point
\[
(x^*,y_1^*,\ldots,y_P^*,\tau^*)
\]
is feasible for the corresponding
program \eqref{lp}. Therefore,
\[
\min_{\boldsymbol r\in\mathcal{R}}
\tau_{\boldsymbol r}
\leq\tau_{\boldsymbol r^*}
\leq\tau^*.
\]
Combining this inequality with
\eqref{eq:lp-first-inequality} establishes
\eqref{eq:finite-lp-equivalence}.

\item
Let $\boldsymbol r^*$ attain the minimum
in \eqref{eq:finite-lp-equivalence}, and
let
$(\bar x,\bar y_1,\ldots,\bar y_P,\bar t)$
be an optimal solution of its linear
program. Then
\[
\bar t=\tau_{\boldsymbol r^*}=\tau^*.
\]
By feasibility,
\[
\Phi(\bar x)\leq\bar t=\tau^*.
\]
Since $\bar x\in S$ and $\tau^*$ is the
minimum of $\Phi$ over $S$, the reverse
inequality also holds. Hence
$\Phi(\bar x)=\tau^*$, and therefore
$\bar x\in X^*$.
\end{enumerate}
\end{proof}
%%%%%%%%%%%%%%%%%%%%%%%%%%%%%%%%%%%%%%%%%%%%%%%%%%%%%%%%%%%%%%%%%%%%%%%%%%%%%%%%%%%%%%%%%%%%%%%%%%%%%%%%%%%%%%%%%%%%%%%%%%%%%%%%%%%%%%%%%%%%%%%%%%%%%%%%%%%%%%%%%%%%%%%%%%%%%%%%%%%%%%%%%%%%%%%%%%%%%%%%%%%%%%%%%%%%%%%%%%%%%%%%%%%%%%%%%%%%%%%%%%%%%%%%%%%%%%%%%%%%%%%%%%%%%%%%%%%%%%%%%%%%%%%%%%%%%%%%%%%%%

\subsection{An Interactive Compromise Method Based on Efficient-Face Generation}
\label{subsec:solution-procedure}

The construction of the efficient sets $E_p$ may require
the examination of a large number of faces of the feasible
polytope. To avoid completing this construction before
computing a compromise, we combine Say{\i}n's facial
decomposition method with the proposed minimax formulation.
The decomposition is interrupted whenever one new efficient
face has been identified for each decision maker. A compromise
is then computed and tested against prescribed tolerances.
If the test fails, the decomposition is resumed.

For each decision maker
$p=1,\ldots,P$, consider the linear objective functions
\[
f_{pj}(x)=c_{pj}^{\top}x,
\qquad j=1,\ldots,q_p,
\]
which are assumed to be minimized.

\medskip
\noindent
\textbf{Efficient-face generation.}
For an index set $I\subseteq\{1,\ldots,m\}$, define
\[
F(I)=\{x\in S:A_Ix=b_I\},
\qquad F(\varnothing)=S.
\]

Following the top-down facial decomposition approach of
Say{\i}n, candidate faces are examined through their defining
active constraints. Their efficiency can be checked by solving
\begin{equation}
\label{eq:face-efficiency-test}
\begin{aligned}
\eta_p(I)=
\max_{u,v}\quad&
\sum_{j=1}^{q_p}
\bigl(f_{pj}(v)-f_{pj}(u)\bigr)\\
\text{subject to}\quad&
u\in S,\quad v\in F(I),\\
&f_{pj}(u)\leq f_{pj}(v),
\qquad j=1,\ldots,q_p.
\end{aligned}
\end{equation}

If \eqref{eq:face-efficiency-test} is infeasible,
then $F(I)=\varnothing$. Otherwise, $F(I)$ is efficient
if and only if $\eta_p(I)=0$. When $\eta_p(I)>0$,
the search continues among its subfaces.

At iteration $k$, the procedure provides one previously
unselected efficient face $Q_p^k$ for each decision maker.
Faces are explored from higher to potentially lower
dimensions, without requiring the complete construction
of $E_p$.

\medskip
\noindent
\textbf{Compromise computation.}
For the selected faces
$Q_1^k,\ldots,Q_P^k$, we solve
\begin{equation}
\label{eq:progressive-compromise}
\begin{aligned}
\tau_k=\min_{x\in S}\;
\max_{1\leq p\leq P}
\alpha_p d_\sigma(x,Q_p^k).
\end{aligned}
\end{equation}

Using the weighted maximum norm, this problem is
equivalent to the linear program
\begin{equation}
\label{eq:progressive-compromise-lp}
\begin{aligned}
\min_{x,y_1,\ldots,y_P,t}\quad &t\\
\text{subject to}\quad
&x\in S,\\
&y_p\in Q_p^k,
&&p=1,\ldots,P,\\
&-\frac{\sigma_i}{\alpha_p}t
\leq x_i-y_{pi}
\leq\frac{\sigma_i}{\alpha_p}t,
&&\substack{i=1,\ldots,n,\\p=1,\ldots,P},\\
&t\geq0.
\end{aligned}
\end{equation}

Let $X_k$ be the set of optimal decision vectors of
\eqref{eq:progressive-compromise-lp}. After $k$
iterations, the best radius and its associated
compromises are given by
\begin{equation}
\label{eq:progressive-incumbent}
U_k=\min_{1\leq r\leq k}\tau_r,
\qquad
\mathcal C_k=
\bigcup_{\substack{1\leq r\leq k\\\tau_r=U_k}}X_r.
\end{equation}

\medskip
\noindent
\textbf{Acceptance criterion and secondary selection.}
For each objective, let
\[
z_{pj}^{\mathrm{id}}
=\min_{x\in S}f_{pj}(x)
\]
be its ideal value. The decision makers specify
nonnegative tolerances $\varepsilon_{pj}$ representing
the largest acceptable deviations from these values.

The set of acceptable decisions is defined by
\begin{equation}
\label{eq:acceptable-set}
\mathcal A_{\varepsilon}
=
\left\{
x\in S:
f_{pj}(x)-z_{pj}^{\mathrm{id}}
\leq\varepsilon_{pj},
\ \forall p,j
\right\}.
\end{equation}

A compromise is accepted if it belongs to
$\mathcal A_{\varepsilon}$. When several acceptable
compromises attain the current best radius, a secondary
linear criterion $g(x)=c^\top x$ is used to select
\begin{equation}
\label{eq:progressive-secondary}
\bar x^k\in
\operatorname*{arg\,min}_{x\in
\mathcal C_k\cap\mathcal A_{\varepsilon}}
g(x).
\end{equation}

This selection amounts to solving linear programs over
the relevant optimal faces of
\eqref{eq:progressive-compromise-lp}, with $t=U_k$
and the additional constraints defining
$\mathcal A_{\varepsilon}$.

The resulting procedure is stated below.

\medskip
\noindent
\textbf{Algorithm 1. Progressive facial decomposition
and compromise selection}

\smallskip
\noindent
\textit{Input:}
The feasible polytope $S$, the objective functions
$f_{pj}$, weights $\alpha_p>0$, scales $\sigma_i>0$,
tolerances $\varepsilon_{pj}\geq0$, and an optional
secondary linear criterion $g$.

\begin{enumerate}

\item \textbf{Initialization.}

Compute the ideal values $z_{pj}^{\mathrm{id}}$
and define $\mathcal A_{\varepsilon}$ by
\eqref{eq:acceptable-set}.

If $\mathcal A_{\varepsilon}=\varnothing$,
terminate: the prescribed tolerances are incompatible.

For each decision maker $p$, initialize
the facial decomposition lists
\[
\mathcal L_p=\{\varnothing\},
\qquad
\mathcal G_p=\varnothing,
\qquad
\mathcal N_p=\varnothing,
\]
where $\mathcal L_p$ contains candidate index sets,
$\mathcal G_p$ contains the identified efficient faces,
and $\mathcal N_p$ records infeasible faces.

Set $k=0$ and $U_0=+\infty$.

\item \textbf{Generation of efficient faces.}

For each $p=1,\ldots,P$, apply the following
steps until one new efficient face is found:

\begin{enumerate}
\item Select an index set $I\in\mathcal L_p$
of minimum cardinality and remove it from
$\mathcal L_p$.

\item Discard $I$ if it contains an index set
already recorded as efficient or infeasible.

\item Otherwise, solve
\eqref{eq:face-efficiency-test}.

\begin{itemize}
\item If the problem is infeasible,
add $I$ to $\mathcal N_p$.

\item If $\eta_p(I)=0$, record $F(I)$
as efficient and set $Q_p^{k+1}=F(I)$.

\item If $\eta_p(I)>0$, insert the
previously unexamined index sets
$I\cup\{i\}$, $i\notin I$, into
$\mathcal L_p$.
\end{itemize}
\end{enumerate}

Once $Q_p^{k+1}$ has been obtained,
suspend the decomposition for decision maker $p$.

If no new efficient face can be generated
for some decision maker, proceed to Step 6.

\item \textbf{Compromise computation.}

Set $k\leftarrow k+1$ and solve
\eqref{eq:progressive-compromise-lp}
using the newly generated faces
$Q_1^k,\ldots,Q_P^k$.

Obtain $\tau_k$ and $X_k$.

\item \textbf{Update.}

Update the best radius and candidate set:
\[
U_k=\min\{U_{k-1},\tau_k\},
\]
\[
\mathcal C_k=
\begin{cases}
X_k,
&\tau_k<U_{k-1},\\[1mm]
\mathcal C_{k-1}\cup X_k,
&\tau_k=U_{k-1},\\[1mm]
\mathcal C_{k-1},
&\tau_k>U_{k-1}.
\end{cases}
\]

\item \textbf{Acceptance test and selection.}

Determine whether
\[
\mathcal C_k\cap\mathcal A_{\varepsilon}
\neq\varnothing.
\]

If this condition holds, select
$\bar x^k$ by \eqref{eq:progressive-secondary},
return $(\bar x^k,U_k)$, and stop.

Otherwise, resume the suspended facial
decomposition procedures and return to Step 2.

\item \textbf{Termination without acceptance.}

If no further face tuple can be generated,
return $U_k$ and $\mathcal C_k$,
and report that the prescribed acceptance
criterion has not been satisfied.

\end{enumerate}

\begin{proposition}
\label{prop:progressive-validity}
Suppose that $S$ is a nonempty compact polytope
and that every generated face $Q_p^k$ is
nonempty and efficient for decision maker $p$.
Then each restricted compromise problem admits
an optimal solution. Moreover,
\[
\tau^*\leq U_{k+1}\leq U_k
\]
at every completed iteration.

If Algorithm 1 terminates at Step 5,
the returned decision is feasible and
satisfies all prescribed objective tolerances.
\end{proposition}

\begin{proof}
Each $Q_p^k$ is a nonempty compact subset
of $S$. Hence the objective function in
\eqref{eq:progressive-compromise} is continuous
on a compact set and attains its minimum.

Since $Q_p^k\subseteq E_p$, we have
\[
d_\sigma(x,E_p)
\leq d_\sigma(x,Q_p^k),
\qquad x\in S.
\]
Consequently, $\tau^*\leq\tau_k$ for every $k$.
Taking the minimum over the iterations already
performed yields
\[
\tau^*\leq U_{k+1}\leq U_k.
\]
Finally, termination at Step 5 ensures that
$\bar x^k\in\mathcal C_k\cap
\mathcal A_{\varepsilon}$, which proves the
last assertion.
\end{proof}

\begin{remark}
    The procedure can terminate before the complete
efficient sets have been constructed. In this case,
the returned decision satisfies the specified
acceptance criterion but is not necessarily a
global minimizer of the original compromise problem.
Indeed, the successive face tuples do not generally
cover all possible combinations of efficient faces.
\end{remark}

%%%%%%%%%%%%%%%%%%%%%%%%%%%%%%%%%%%%%%%%%%%%%%%%%%%%%%%%%%%%%%%%%%%%%%ùùùùùùùùùùùùùùùùùùùùùùùùùùùùùùùùùùùùùùùùù%%%%%%%%%%%%%%%%%%%%%%%%%%%%%%%%%%%%%%%%%%%%%%%%%%%%%%%%%%%%%%%%%%%%%%%%%%%%%%%%%%%%%%%%%%%%%%%%%%%%%%%%%%%%%%%%%%%%%%%%%%%%%%%%%%%%%%%%%%%%%%%%%%%%%%%%%%%%%%%%%%%%%%%%%%%%%%%%%%%%%%%%%%%%%%%%%%%%%%%%%%%%%%%%%%%%%%%%%%%%%%%%%%%

\subsection{Stability under set approximation}\label{secapp}
Approximate efficient sets can be used when a verified geometric error bound is available. Write $d_{H,\sigma}$ for the Hausdorff distance induced by $\|\cdot\|_\sigma$.
\begin{theorem}[Approximation error]
Let $\widetilde E_p$ be nonempty compact sets satisfying
$d_{H,\sigma}(E_p,\widetilde E_p)\le\eta_p$. On the same domain $S$, define
\[
\widetilde\Phi(x)=\max_p\alpha_p\dist_\sigma(x,\widetilde E_p),\quad
\widetilde\tau^*=\min_S\widetilde\Phi,\quad
\varepsilon=\max_p\alpha_p\eta_p.
\]
Then $|\widetilde\tau^*-\tau^*|\le\varepsilon$. If $\widetilde x$ minimizes $\widetilde\Phi$ on $S$, then
$\Phi(\widetilde x)\le\tau^*+2\varepsilon$.
\end{theorem}
\begin{proof}
The Hausdorff bound and triangle inequality imply
$|\dist_\sigma(x,E_p)-\dist_\sigma(x,\widetilde E_p)|\le\eta_p$.
Consequently, $|\Phi(x)-\widetilde\Phi(x)|\le\varepsilon$ uniformly on $S$. Taking minima yields the first assertion. Finally,
\[
\Phi(\widetilde x)\le\widetilde\Phi(\widetilde x)+\varepsilon
=\widetilde\tau^*+\varepsilon\le\tau^*+2\varepsilon.
\]
\end{proof}
This bound concerns values, not distances between minimizers. Small changes in objective coefficients do not automatically imply a small Hausdorff perturbation of the efficient sets.

\section{Numerical example: harvest allocation between two fishing zones}
\label{sec:application}

\subsection{Decision variables and model parameters}

We consider two fishing zones connected by biomass movement over a fixed management period. A common harvest decision must reconcile two concerns: increasing catches in each zone and preserving the corresponding residual biomass.

The decision vector is
\[
u=(u_1,u_2),
\]
where $u_i$ denotes the quantity harvested in zone $i$, expressed in thousand tonnes (kt). Let $B_i$ be the initial biomass in that zone, and let $m_{ij}\in[0,1]$ denote the fraction moving from zone $i$ to zone $j$ before harvest. The available biomasses are therefore
\begin{equation}\label{eq:available_biomass}
\begin{aligned}
A_1&=(1-m_{12})B_1+m_{21}B_2,\\
A_2&=m_{12}B_1+(1-m_{21})B_2.
\end{aligned}
\end{equation}
Movement preserves the aggregate biomass:
\[
A_1+A_2=B_1+B_2=B.
\]
After harvesting, the residual biomass in zone $i$ is
\begin{equation}\label{eq:residual_biomass}
R_i(u)=A_i-u_i,\qquad i=1,2.
\end{equation}
This is a single-period accounting model; growth and natural mortality are not included.

The remaining parameters are the local harvest capacities $\bar u_i$, the minimum residual biomasses $L_i$, and the aggregate harvest ceiling $Q$. All these quantities are fixed before the decision is made.

\subsection{Constraints and objective functions}

A feasible allocation must satisfy three requirements:
\begin{enumerate}
\item Harvest is nonnegative and does not exceed local operational capacity:
\[
0\le u_i\le\bar u_i,\qquad i=1,2.
\]
\item Residual biomass remains above the prescribed local floor:
\[
A_i-u_i\ge L_i,\qquad i=1,2.
\]
\item Total harvest does not exceed the aggregate ceiling:
\[
u_1+u_2\le Q.
\]
\end{enumerate}

Assuming $A_i\ge L_i$, define the effective capacities
\[
c_i=\min\{\bar u_i,A_i-L_i\},\qquad i=1,2.
\]
The common feasible region is the compact polytope
\begin{equation}\label{eq:harvest_feasible}
S=\left\{
u\in\mathbb R^2:
0\le u_i\le c_i,\ i=1,2,\quad u_1+u_2\le Q
\right\}.
\end{equation}
It is nonempty because it contains the zero-harvest decision.

The harvesting authority seeks to maximize catches in both zones:
\begin{equation}\label{eq:harvest_objectives}
\max_{u\in S}F_E(u),
\qquad
F_E(u)=
\begin{pmatrix}
u_1\\u_2
\end{pmatrix}.
\end{equation}
The conservation authority seeks to maximize the two residual biomasses:
\begin{equation}\label{eq:conservation_objectives}
\max_{u\in S}F_C(u),
\qquad
F_C(u)=
\begin{pmatrix}
A_1-u_1\\A_2-u_2
\end{pmatrix}.
\end{equation}
Both vector objectives are maximized componentwise. Since $A_1$ and $A_2$ are constant, the conservation objective is equivalent to maximizing $(-u_1,-u_2)^\top$. The economic criterion measures harvest quantities; no monetary costs or revenues are assigned.

Let
\[
T=\min\{Q,c_1+c_2\}.
\]
The efficient reference sets are
\begin{equation}\label{eq:harvest_reference_sets}
E_E=\{v\in S:v_1+v_2=T\},
\qquad
E_C=\{(0,0)\}.
\end{equation}
Indeed, an allocation with total harvest below $T$ admits an increase in at least one coordinate without decreasing the other. Conversely, no allocation with total $T$ can be dominated by another feasible allocation. For conservation, zero harvest dominates every nonzero feasible decision. Consequently, $E_E\cap E_C=\varnothing$ whenever $T>0$.

\subsection{Compromise model}\label{sec:capped}

Let $\sigma_1,\sigma_2>0$ be fixed distance scales. We use
\[
\|z\|_\sigma
=
\max\left\{
\frac{|z_1|}{\sigma_1},
\frac{|z_2|}{\sigma_2}
\right\}.
\]
The distances from a feasible decision to the two reference sets are
\[
d_E(u)=\min_{v\in E_E}\|u-v\|_\sigma,
\qquad
d_C(u)=\|u\|_\sigma.
\]
For weights $\alpha,\beta>0$, the compromise problem is
\begin{equation}\label{eq:harvest_compromise}
\min_{u\in S}
\max\{\alpha d_E(u),\beta d_C(u)\}.
\end{equation}
The weight $\alpha$ gives importance to proximity to the harvesting
reference set, while $\beta$ gives importance to proximity
to the conservation reference.

Introducing a reference point $v\in E_E$ gives the linear program
\begin{equation}\label{eq:harvest_linear_program}
\begin{aligned}
\min_{u,v,\tau}\quad &\tau\\
\text{subject to}\quad
&0\le u_i\le c_i,
&&i=1,2,\\
&u_1+u_2\le Q,\\
&0\le v_i\le c_i,
&&i=1,2,\\
&v_1+v_2=T,\\
&-\frac{\sigma_i}{\alpha}\tau
\le u_i-v_i
\le\frac{\sigma_i}{\alpha}\tau,
&&i=1,2,\\
&u_i\le\frac{\sigma_i}{\beta}\tau,
&&i=1,2,\\
&\tau\ge0.
\end{aligned}
\end{equation}
The vector $u$ is the harvest decision.
The vector $v$ is used only to measure its distance
from the harvesting efficient set.

For this allocation model, a closest point of $E_E$ to the origin
can also be computed directly. Let $\delta$ be the smallest
nonnegative number satisfying
\begin{equation}\label{eq:harvest_threshold}
\min\{c_1,\delta\sigma_1\}
+
\min\{c_2,\delta\sigma_2\}
=T.
\end{equation}
Then
\begin{equation}\label{eq:harvest_selection}
v_i^\circ=\min\{c_i,\delta\sigma_i\},
\qquad
u^\circ=\frac{\alpha}{\alpha+\beta}v^\circ,
\qquad
\tau^*=\frac{\alpha\beta}{\alpha+\beta}\delta.
\end{equation}
To see why $\delta$ is the distance between the reference sets,
observe that $\|v\|_\sigma\le t$ requires
$v_i\le\min\{c_i,t\sigma_i\}$.
Such a point can have total harvest $T$ precisely when the sum
of these bounds is at least $T$.

Formula~\eqref{eq:harvest_selection} then follows from the
two-set compromise construction. We use this formula to specify
the reported decision when several optimal compromises exist.
Its optimal radius is checked independently by solving
\eqref{eq:harvest_linear_program}.
\subsection{Numerical data}

The numerical scale is taken from the 2025 ICES assessment of northern Norwegian coastal cod~\cite{ices}. The reported biomass and catch for 2024, converted to kt, are
\begin{equation}\label{eq:harvest_data}
B=276.047506343027,
\qquad
Q=28.487.
\end{equation}
The reported catch is used here as a scenario ceiling. The spatial split, movement fractions, capacities and retention floors are prescribed model assumptions.

\begin{table}[htbp]
\centering
\caption{Parameters of the numerical scenario. Biomass and harvest quantities are in kt.}
\label{tab:harvest_parameters}
\begin{tabular}{lll}
\hline
Parameter & Specification & Meaning\\
\hline
$B_1,B_2$ & $0.6B,\;0.4B$ & Initial biomass\\
$m_{12},m_{21}$ & $0.10,\;0.20$ & Movement fractions\\
$\bar u_1,\bar u_2$ & $0.6Q,\;0.8Q$ & Operational capacities\\
$L_1,L_2$ & $0.75A_1,\;0.75A_2$ & Residual biomass floors\\
$\sigma_1,\sigma_2$ & $\bar u_1,\;\bar u_2$ & Distance scales\\
$\alpha,\beta$ & $1,\;1$ & Authority weights\\
\hline
\end{tabular}
\end{table}

Equations~\eqref{eq:available_biomass} give
\[
A_1=0.62B,\qquad A_2=0.38B.
\]
The resulting local quantities are shown in
Table~\ref{tab:harvest_accounting}.

\begin{table}[htbp]
\centering
\caption{Local biomass accounting and effective harvest capacities, in kt. Values are rounded to six decimal places.}
\label{tab:harvest_accounting}
\begin{tabular}{lrr}
\hline
Quantity & Zone 1 & Zone 2\\
\hline
Initial biomass $B_i$ & 165.628504 & 110.419003\\
Available biomass $A_i$ & 171.149454 & 104.898052\\
Residual floor $L_i$ & 128.362090 & 78.673539\\
Removable biomass $A_i-L_i$ & 42.787363 & 26.224513\\
Operational capacity $\bar u_i$ & 17.092200 & 22.789600\\
Effective capacity $c_i$ & 17.092200 & 22.789600\\
\hline
\end{tabular}
\end{table}

In both zones, the operational capacity is smaller than the removable biomass. Hence
\[
c_1=0.6Q,\qquad c_2=0.8Q,\qquad T=Q.
\]
The reference sets reduce to
\begin{equation}\label{eq:numerical_reference_sets}
E_E=
\{(s,Q-s):0.2Q\le s\le0.6Q\},
\qquad
E_C=\{(0,0)\}.
\end{equation}

\subsection{Numerical results and sensitivity analysis}\label{sensitivity}

We solved \eqref{eq:harvest_linear_program} using
\texttt{scipy.optimize.linprog} with the HiGHS solver.
All computations used unrounded inputs. In the initial
scenario, with $\alpha=\beta=1$, the solution is
\begin{equation}\label{eq:numerical_compromise}
u^*=
\left(\frac{3Q_0}{14},\frac{2Q_0}{7}\right)
=(6.104357142857,\,8.139142857143),
\qquad
\tau^*=\frac{5}{14}.
\end{equation}
A closest harvesting reference is
\[
v^\circ=
\left(\frac{3Q_0}{7},\frac{4Q_0}{7}\right)
=(12.208714285714,\,16.278285714286).
\]
Indeed, since $\sigma_1+\sigma_2=1.4Q_0$, every feasible
solution satisfies
\[
Q_0=v_1+v_2
\le u_1+u_2+\tau(\sigma_1+\sigma_2)
\le 2\tau(\sigma_1+\sigma_2),
\]
which gives $\tau\ge 5/14$. The proposed points attain this
bound, with
\[
\frac{u_i^*}{\sigma_i}
=\frac{v_i^\circ-u_i^*}{\sigma_i}
=\frac{5}{14},
\qquad i=1,2.
\]
Equality in the summed inequalities forces
$u_i^*=\tau^*\sigma_i$ for both coordinates, proving uniqueness
in this scenario. The total harvest is $14.2435$~kt,
and the remaining biomass is $261.804006$~kt, approximately
$94.84\%$ of the initial total.

\begin{table}[htbp]
\centering\small
\caption{Harvests and remaining biomasses in the initial
scenario. All quantities are in kt.}
\label{tab:harvest_results}
\begin{tabular}{lrrr}
\toprule
Quantity & Conservation & Compromise & Harvest reference\\
\midrule
$u_1$ & 0.000000 & 6.104357 & 12.208714\\
$u_2$ & 0.000000 & 8.139143 & 16.278286\\
$u_1+u_2$ & 0.000000 & 14.243500 & 28.487000\\
$R_1(u)$ & 171.149454 & 165.045097 & 158.940740\\
$R_2(u)$ & 104.898052 & 96.758910 & 88.619767\\
$R_1(u)+R_2(u)$ & 276.047506 & 261.804006 & 247.560506\\
\bottomrule
\end{tabular}
\end{table}

\begin{figure}[htbp]
\centering
\includegraphics[width=\textwidth]{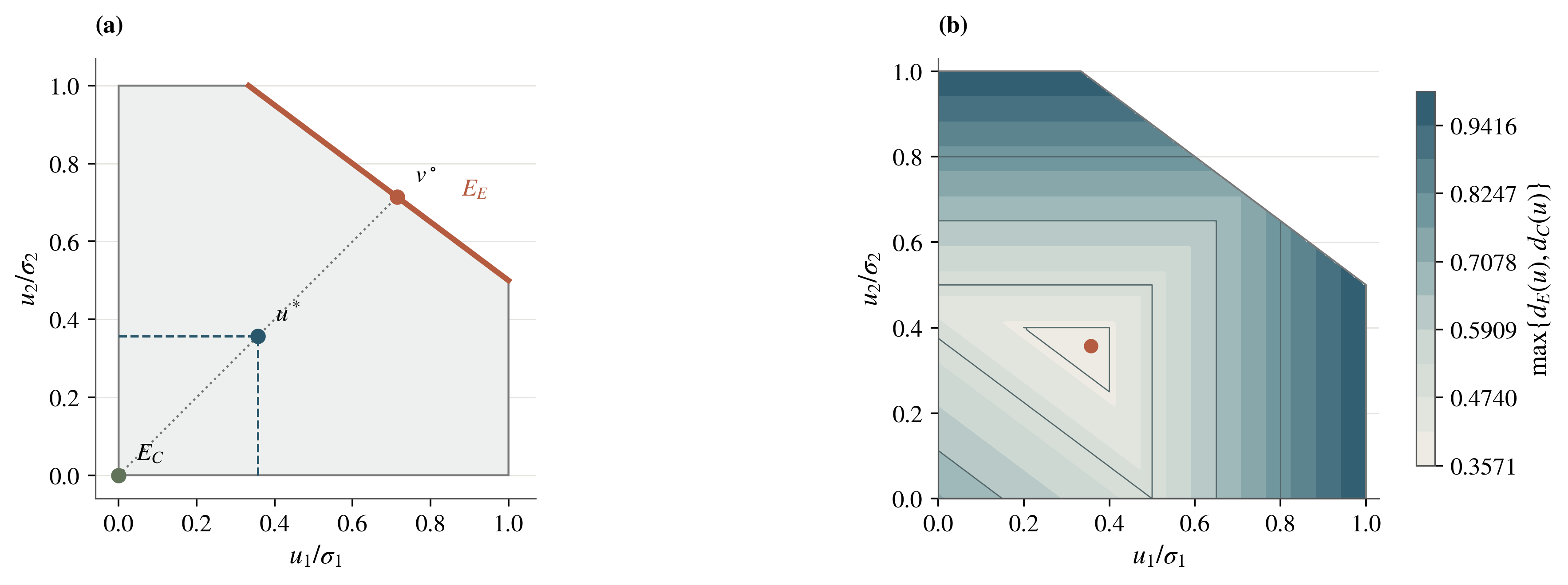}
\caption{Initial scenario in normalized coordinates.
Left: feasible region, harvesting efficient set, closest
reference and compromise. Right: compromise criterion
over the feasible region.}
\label{fig:harvest_geometry}
\end{figure}

We next examine the effects of the authority weights,
migration, distance scales, biomass requirements and total
harvest ceiling. Each experiment varies one parameter while
keeping the remaining inputs fixed as specified below.
For each scenario, let
\[
c_i=\min\{\bar u_i,A_i-L_i\},
\qquad
T=\min\{Q,c_1+c_2\}.
\]
To specify the reported solution when several minimizers exist,
we use
\begin{equation}\label{eq:harvest_selection_sensitivity}
v_i^\circ=\min\{c_i,\delta\sigma_i\},
\qquad
u^\circ=\frac{\alpha}{\alpha+\beta}v^\circ,
\qquad
\tau^*=\frac{\alpha\beta}{\alpha+\beta}\delta,
\end{equation}
where $\delta$ is the smallest nonnegative solution of
\begin{equation}\label{eq:harvest_threshold_sensitivity}
\sum_{i=1}^{2}\min\{c_i,\delta\sigma_i\}=T.
\end{equation}
The threshold is computed by bisection, and the resulting
radius is compared with the independently solved linear program.

\medskip
\noindent\textbf{Authority weights.}
We fix $\beta=1$ and vary $\alpha$ over $[0.1,10]$
using 81 logarithmically spaced values.
The computed solutions agree with
\[
u_{\alpha,1}^*
=\frac{\alpha}{1+\alpha}
\left(\frac{3Q_0}{7},\frac{4Q_0}{7}\right),
\qquad
\tau_{\alpha,1}^*
=\frac{5\alpha}{7(1+\alpha)}.
\]
Thus, increasing the harvesting weight raises the total
harvest from $Q_0/11$ to $10Q_0/11$, while preserving
the allocation ratio $u_1^*/u_2^*=3/4$.

\begin{table}[htbp]
\centering\small
\caption{Effect of the harvesting weight, with $\beta=1$.
Physical quantities are in kt.}
\label{tab:weight_variation}
\begin{tabular}{lrrrrrr}
\toprule
$\alpha$ & $c_1$ & $c_2$ & $u_1^\circ$ & $u_2^\circ$ & $T$ & $\tau^*$ \\
\midrule
0.10 & 17.092200 & 22.789600 & 1.109883 & 1.479844 & 28.487000 & 0.064935 \\
0.25 & 17.092200 & 22.789600 & 2.441743 & 3.255657 & 28.487000 & 0.142857 \\
0.50 & 17.092200 & 22.789600 & 4.069571 & 5.426095 & 28.487000 & 0.238095 \\
1.00 & 17.092200 & 22.789600 & 6.104357 & 8.139143 & 28.487000 & 0.357143 \\
2.00 & 17.092200 & 22.789600 & 8.139143 & 10.852190 & 28.487000 & 0.476190 \\
4.00 & 17.092200 & 22.789600 & 9.766971 & 13.022629 & 28.487000 & 0.571429 \\
10.00 & 17.092200 & 22.789600 & 11.098831 & 14.798442 & 28.487000 & 0.649351 \\
\bottomrule
\end{tabular}

\end{table}

\begin{figure}[htbp]
\centering
\includegraphics[width=\textwidth]{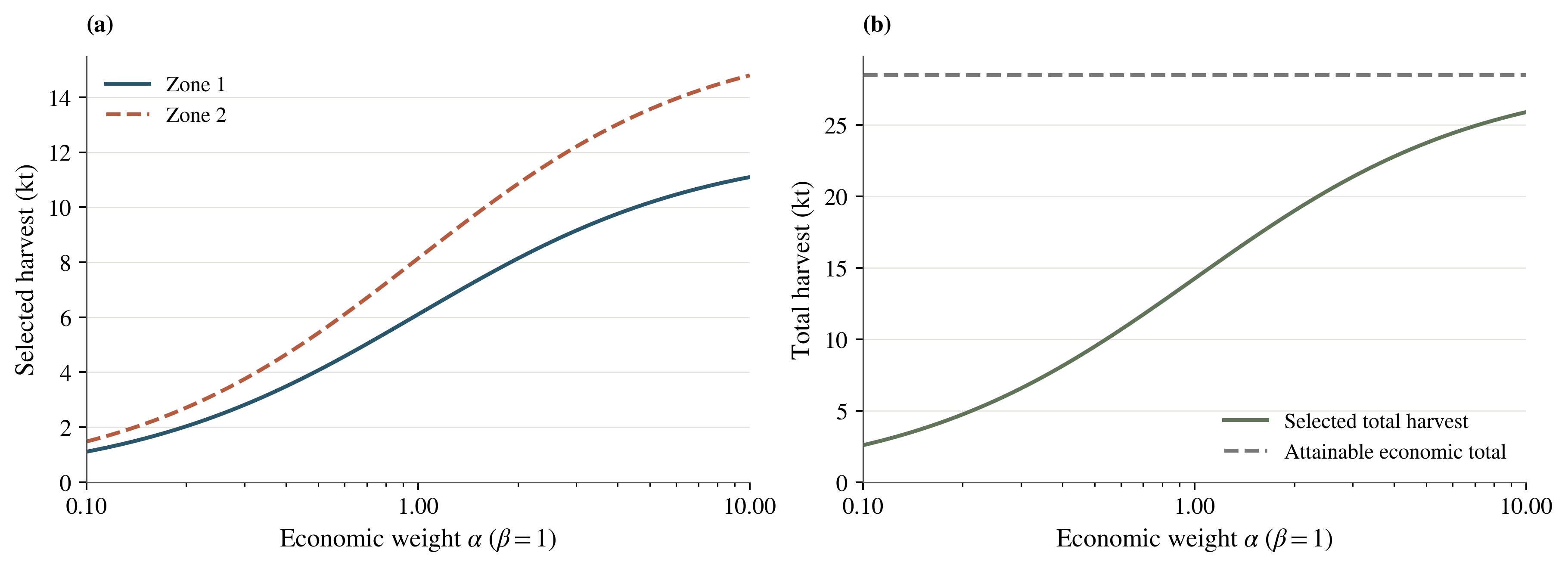}
\caption{Local and total harvests as functions of $\alpha$.
The horizontal axis is logarithmic.}
\label{fig:weight_variation}
\end{figure}

\medskip
\noindent\textbf{Migration.}
To examine situations where biomass availability affects
the decision, we impose the fixed requirements
$L_i=0.85A_i^0$, where
$A_1^0=0.62B$ and $A_2^0=0.38B$.
We retain $m_{21}=0.20$ and vary $m=m_{12}$ over
$[0.03,0.25]$ using 111 equally spaced values. Then
\[
A_1(m)=B(0.68-0.6m),
\qquad
A_2(m)=B(0.32+0.6m).
\]
All tested scenarios are feasible. The weights remain equal
and the scales remain $(0.6Q_0,0.8Q_0)$.
Low migration limits the biomass available for harvesting
in zone~2, whereas high migration makes zone~1 restrictive.
The selected total satisfies
\[
u_1^\circ(m)+u_2^\circ(m)
=\frac12\min\{Q_0,c_1(m)+c_2(m)\}.
\]

\begin{table}[htbp]
\centering\small
\caption{Effect of migration with fixed requirements
$L_i=0.85A_i^0$. Physical quantities are in kt.}
\label{tab:migration_variation}
\begin{tabular}{lrrrrrr}
\toprule
$m_{12}$ & $c_1$ & $c_2$ & $u_1^\circ$ & $u_2^\circ$ & $T$ & $\tau^*$ \\
\midrule
0.03 & 17.092200 & 4.140713 & 8.546100 & 2.070356 & 21.232913 & 0.500000 \\
0.06 & 17.092200 & 9.109568 & 8.546100 & 4.554784 & 26.201768 & 0.500000 \\
0.10 & 17.092200 & 15.734708 & 6.376146 & 7.867354 & 28.487000 & 0.373044 \\
0.15 & 17.092200 & 22.789600 & 6.104357 & 8.139143 & 28.487000 & 0.357143 \\
0.20 & 9.109568 & 22.789600 & 4.554784 & 9.688716 & 28.487000 & 0.425138 \\
0.25 & 0.828143 & 22.789600 & 0.414071 & 11.394800 & 23.617743 & 0.500000 \\
\bottomrule
\end{tabular}

\end{table}

\begin{figure}[htbp]
\centering
\includegraphics[width=\textwidth]{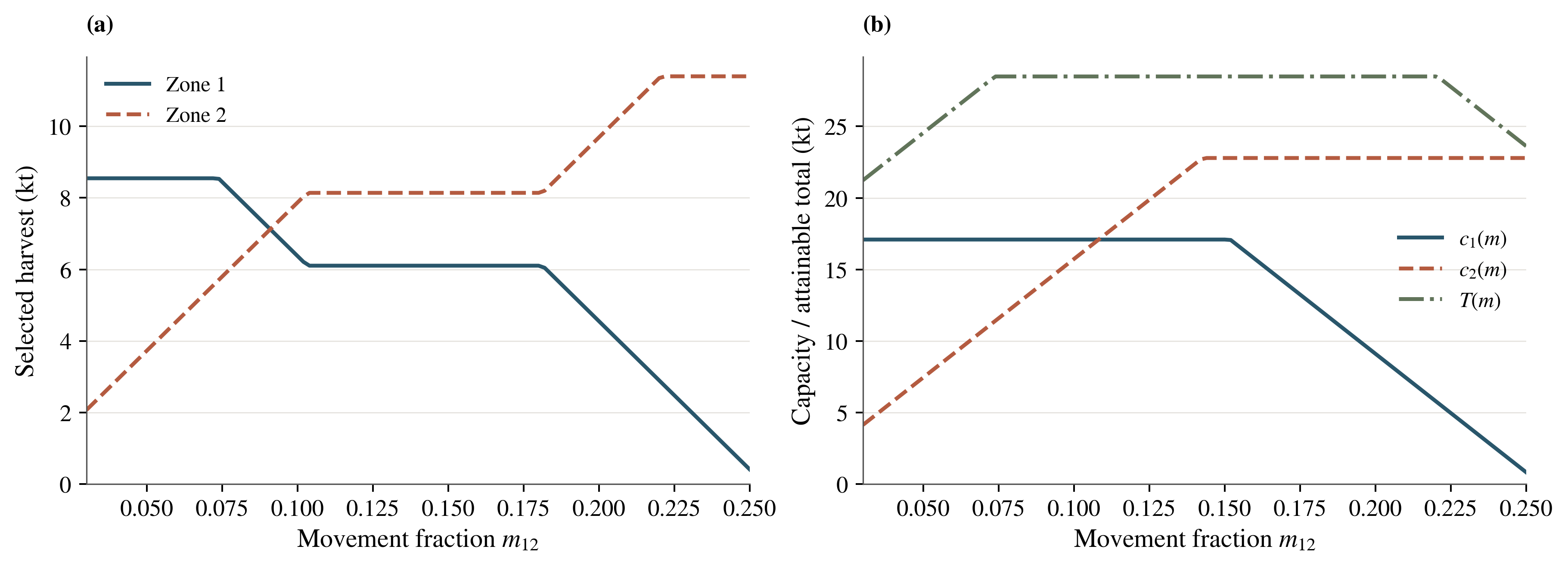}
\caption{Effect of migration on local harvests, effective
capacities and the largest feasible total harvest.}
\label{fig:migration_variation}
\end{figure}

\medskip
\noindent\textbf{Distance scales.}
Under the initial physical constraints and equal weights,
we fix $\sigma_2=Q_0$ and vary $r=\sigma_1/\sigma_2$
over $[0.1,4]$ using 81 logarithmically spaced values.
The selected nearest reference is
\[
v^\circ(r)=
\begin{cases}
(0.2Q_0,0.8Q_0), & 0<r\le 1/4,\\[1mm]
\left(\dfrac{rQ_0}{1+r},\dfrac{Q_0}{1+r}\right),
& 1/4<r<3/2,\\[2mm]
(0.6Q_0,0.4Q_0), & r\ge 3/2.
\end{cases}
\]
Since $u^\circ=v^\circ/2$, the total remains $Q_0/2$.
Increasing $r$ shifts harvest toward zone~1 until the
reference reaches its local capacity. Each radius is
evaluated in the norm associated with the specified scales.

\begin{table}[htbp]
\centering\small
\caption{Effect of the distance scales, with $\sigma_2=Q_0$.
Physical quantities are in kt.}
\label{tab:scale_variation}
\begin{tabular}{lrrrrrr}
\toprule
$\sigma_1/\sigma_2$ & $c_1$ & $c_2$ & $u_1^\circ$ & $u_2^\circ$ & $T$ & $\tau^*$ \\
\midrule
0.10 & 17.092200 & 22.789600 & 2.848700 & 11.394800 & 28.487000 & 1.000000 \\
0.25 & 17.092200 & 22.789600 & 2.848700 & 11.394800 & 28.487000 & 0.400000 \\
0.50 & 17.092200 & 22.789600 & 4.747833 & 9.495667 & 28.487000 & 0.333333 \\
0.75 & 17.092200 & 22.789600 & 6.104357 & 8.139143 & 28.487000 & 0.285714 \\
1.00 & 17.092200 & 22.789600 & 7.121750 & 7.121750 & 28.487000 & 0.250000 \\
2.00 & 17.092200 & 22.789600 & 8.546100 & 5.697400 & 28.487000 & 0.200000 \\
4.00 & 17.092200 & 22.789600 & 8.546100 & 5.697400 & 28.487000 & 0.200000 \\
\bottomrule
\end{tabular}

\end{table}

\begin{figure}[htbp]
\centering
\includegraphics[width=\textwidth]{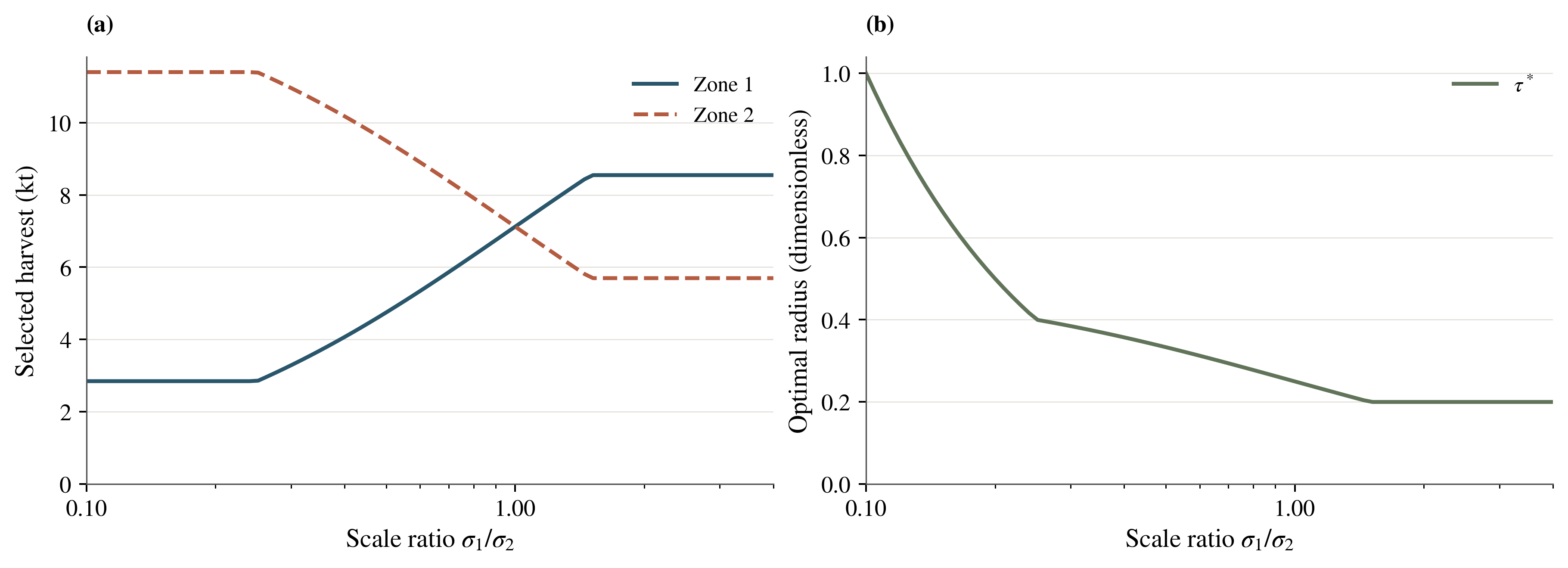}
\caption{Selected local harvests and optimal radius
as functions of the scale ratio.}
\label{fig:scale_variation}
\end{figure}

\medskip
\noindent\textbf{Biomass requirements and harvest ceiling.}
We first set $L_i=\rho A_i^0$ and vary $\rho$ over
$[0.75,0.99]$ using 97 equally spaced values.
The effective capacities become
\[
c_i(\rho)=\min\{\bar u_i,(1-\rho)A_i^0\}.
\]
The compromise initially remains unchanged. As the requirements
become restrictive, the allocation changes, and the selected
total decreases once $c_1(\rho)+c_2(\rho)<Q_0$.
A local harvest may initially increase to compensate for a
reduction in the other zone, before both harvests decrease
under stronger requirements.

\begin{table}[htbp]
\centering\small
\caption{Effect of the required remaining biomass fraction.
Physical quantities are in kt.}
\label{tab:retention_variation}
\begin{tabular}{lrrrrrr}
\toprule
$\rho$ & $c_1$ & $c_2$ & $u_1^\circ$ & $u_2^\circ$ & $T$ & $\tau^*$ \\
\midrule
0.75 & 17.092200 & 22.789600 & 6.104357 & 8.139143 & 28.487000 & 0.357143 \\
0.80 & 17.092200 & 20.979610 & 6.104357 & 8.139143 & 28.487000 & 0.357143 \\
0.85 & 17.092200 & 15.734708 & 6.376146 & 7.867354 & 28.487000 & 0.373044 \\
0.90 & 17.092200 & 10.489805 & 8.546100 & 5.244903 & 27.582005 & 0.500000 \\
0.95 & 8.557473 & 5.244903 & 4.278736 & 2.622451 & 13.802375 & 0.250333 \\
0.99 & 1.711495 & 1.048981 & 0.855747 & 0.524490 & 2.760475 & 0.050067 \\
\bottomrule
\end{tabular}

\end{table}

\begin{figure}[htbp]
\centering
\includegraphics[width=\textwidth]{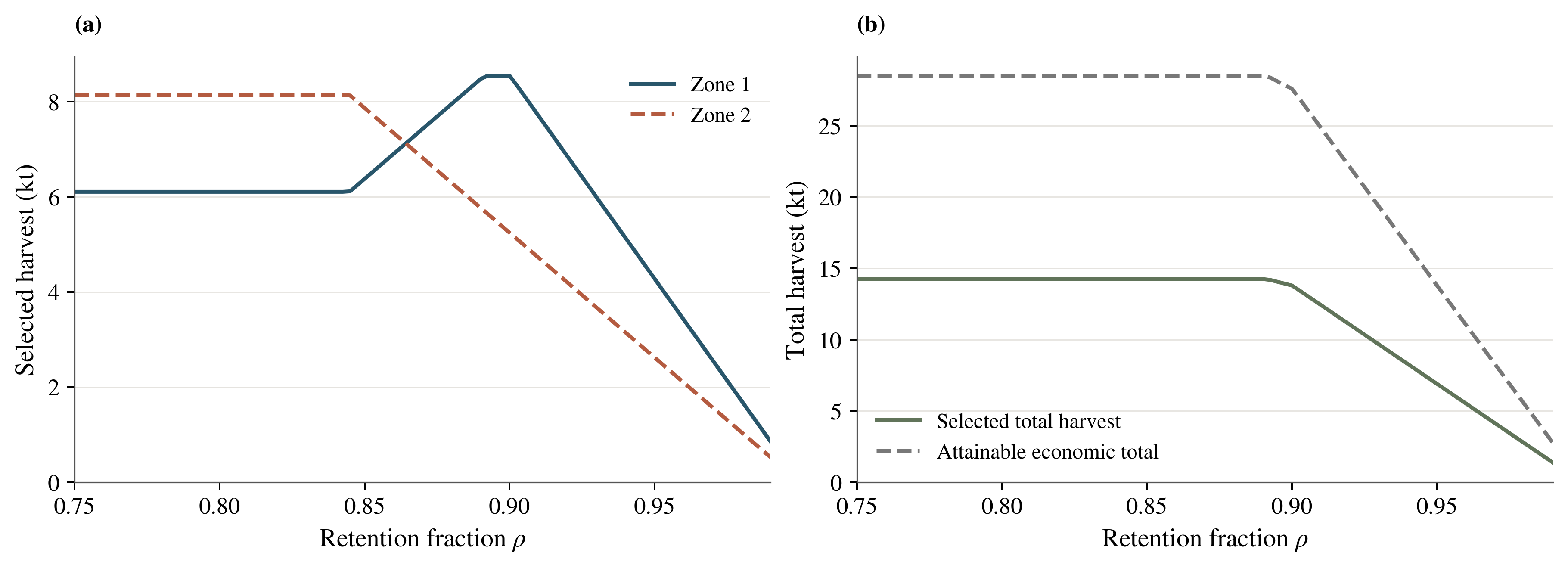}
\caption{Effect of biomass requirements on local
and total harvests.}
\label{fig:retention_variation}
\end{figure}

Finally, we vary $Q=qQ_0$ over $q\in[0.2,1.6]$
using 85 equally spaced values, with the initial biomass
requirements and fixed capacities
$(\bar u_1,\bar u_2)=(0.6Q_0,0.8Q_0)$.
The weights and scales retain their initial values.
Consequently,
\[
T(q)=Q_0\min\{q,1.4\},
\qquad
u_1^\circ+u_2^\circ=\frac{Q_0}{2}\min\{q,1.4\}.
\]
The selected total increases with the ceiling until
$q=1.4$, when both local capacities determine the largest
feasible harvest.

\begin{table}[htbp]
\centering\small
\caption{Effect of the total harvest ceiling with fixed
local capacities. Physical quantities are in kt.}
\label{tab:ceiling_variation}
\begin{tabular}{lrrrrrr}
\toprule
$Q/Q_0$ & $c_1$ & $c_2$ & $u_1^\circ$ & $u_2^\circ$ & $T$ & $\tau^*$ \\
\midrule
0.20 & 17.092200 & 22.789600 & 1.220871 & 1.627829 & 5.697400 & 0.071429 \\
0.50 & 17.092200 & 22.789600 & 3.052179 & 4.069571 & 14.243500 & 0.178571 \\
0.75 & 17.092200 & 22.789600 & 4.578268 & 6.104357 & 21.365250 & 0.267857 \\
1.00 & 17.092200 & 22.789600 & 6.104357 & 8.139143 & 28.487000 & 0.357143 \\
1.20 & 17.092200 & 22.789600 & 7.325229 & 9.766971 & 34.184400 & 0.428571 \\
1.40 & 17.092200 & 22.789600 & 8.546100 & 11.394800 & 39.881800 & 0.500000 \\
1.60 & 17.092200 & 22.789600 & 8.546100 & 11.394800 & 39.881800 & 0.500000 \\
\bottomrule
\end{tabular}

\end{table}

\begin{figure}[htbp]
\centering
\includegraphics[width=\textwidth]{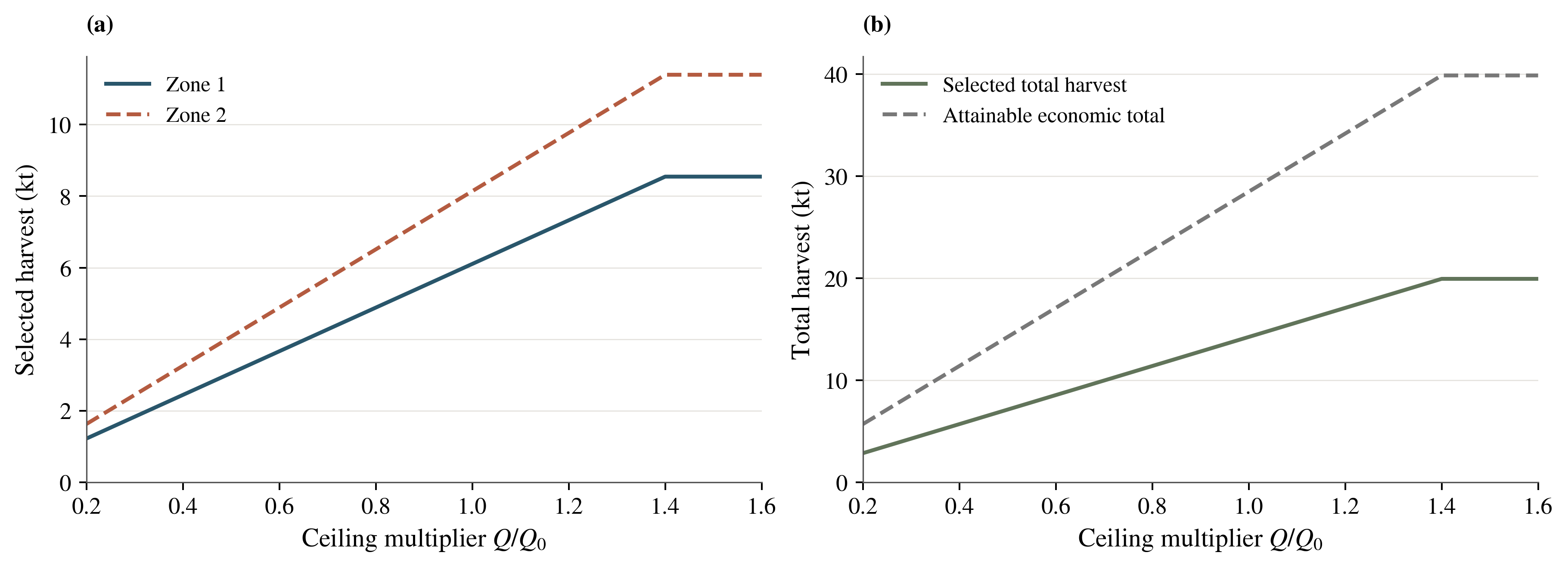}
\caption{Local and total harvests as functions
of the aggregate ceiling.}
\label{fig:ceiling_variation}
\end{figure}

\subsection{Computational verification and discussion}

The five experiments comprise 455 parameter values.
Each linear program has five variables, seven inequality
constraints, one equality constraint and variable bounds.
Calculations used Python~3.12.14, NumPy~2.3.5 and
SciPy~1.17.0 in a Linux environment with an 8~GiB
memory limit. HiGHS and the numerical libraries used
one thread. Each experiment was repeated in five independent
Python processes, with three preliminary solves per process.
Solver calls were timed using \texttt{perf\_counter\_ns},
excluding imports, matrix preparation, analytical checks,
plotting and file output.

\begin{table}[htbp]
\centering\small
\caption{Execution times over five repetitions.
$N$ denotes the number of parameter values.
The final column gives the median total solve time
for one complete experiment.}
\label{tab:execution_times}
\begin{tabular}{lrrrrr}
\toprule
Sweep & $N$ & Median (ms) & IQR (ms) & 95th (ms) & Batch (ms) \\
\midrule
Weights & 81 & 1.112 & 0.989--1.255 & 1.633 & 96.44 \\
Migration & 111 & 1.089 & 0.976--1.250 & 1.822 & 128.66 \\
Scales & 81 & 1.046 & 0.959--1.158 & 1.464 & 89.90 \\
Retention & 97 & 1.147 & 1.022--1.322 & 1.690 & 115.32 \\
Ceiling & 85 & 1.094 & 1.016--1.220 & 1.463 & 97.37 \\
\bottomrule
\end{tabular}

\end{table}

Memory was measured through \texttt{/proc/self/status}.
The measurements cover Python, the numerical libraries,
the solver and the prepared matrices, with the plotting
library excluded. Table~\ref{tab:memory_use} reports
resident memory after the preliminary solves, the largest
recorded peak and its largest increase during timed calculations.

\begin{table}[htbp]
\centering\small
\caption{Process memory over five repetitions, in MiB
($1~\mathrm{MiB}=2^{20}$ bytes).}
\label{tab:memory_use}
\begin{tabular}{lrrr}
\toprule
Sweep & Warm RSS (MiB) & Peak RSS (MiB) & Peak rise (MiB) \\
\midrule
Weights & 59.57 & 59.64 & 0.06 \\
Migration & 59.66 & 59.73 & 0.07 \\
Scales & 59.57 & 59.64 & 0.06 \\
Retention & 59.64 & 59.70 & 0.06 \\
Ceiling & 59.59 & 59.65 & 0.05 \\
\bottomrule
\end{tabular}

\end{table}

\begin{figure}[htbp]
\centering
\includegraphics[width=\textwidth]{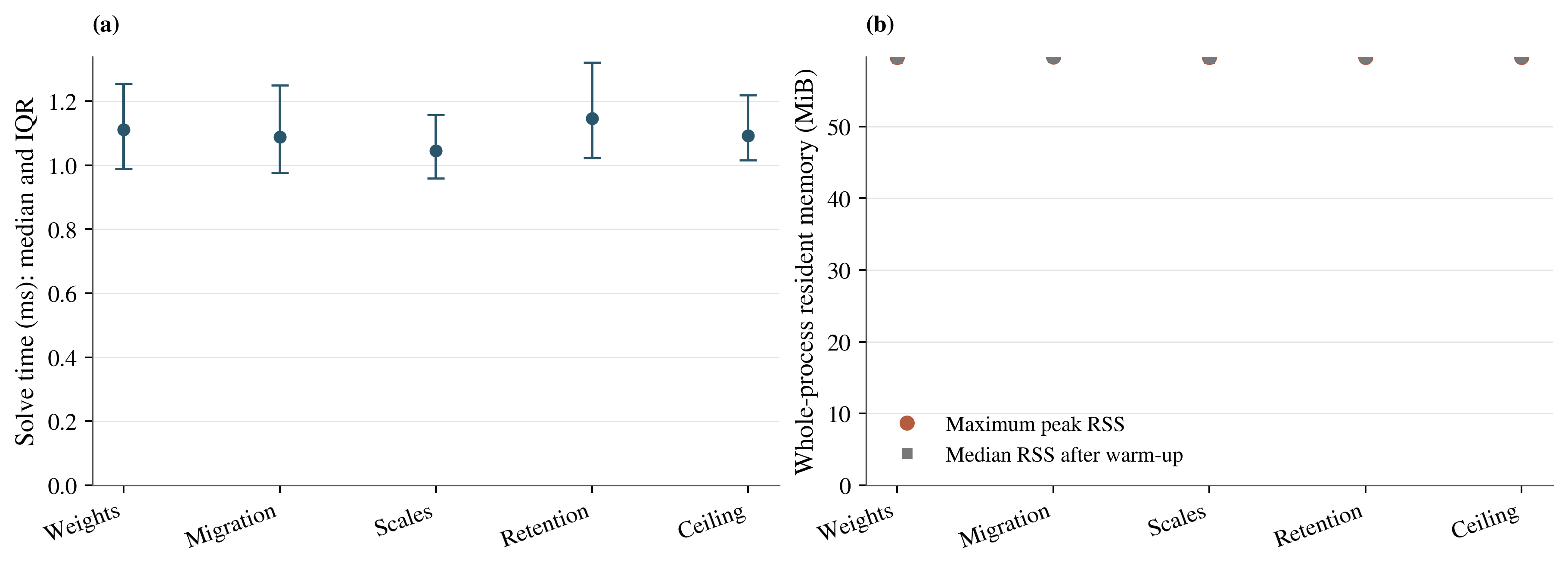}
\caption{Computational measurements.
Left: median solve times and interquartile intervals.
Right: maximum memory peaks and median memory
after preliminary solves.}
\label{fig:computational_results}
\end{figure}

Across all 455 scenarios, the largest absolute difference
between the linear programming radius and
\eqref{eq:harvest_selection} was $3.33\times10^{-16}$.
The largest constraint violation for the selected decisions
and reference points was $1.42\times10^{-14}$.
The timing experiment comprised 2275 measured solves
and 75 preliminary solves. Displayed table rows were
computed directly. The accompanying source code, numerical
data and environment record document these calculations.

\medskip
\noindent\textbf{Discussion.}
The experiments clarify how preferences and resource
constraints determine the selected compromise. Increasing
the harvesting weight raises the total catch, while changing
the distance scales modifies its distribution between zones.
Migration influences the decision through the biomass
available above the fixed local requirements. Stronger
requirements can reduce the attainable total, and an increase
in the harvest ceiling becomes ineffective once local
capacities are reached. For equal weights, the selection
\eqref{eq:harvest_selection} satisfies
$u_1^\circ+u_2^\circ=T/2$, which accounts for the totals
observed throughout the experiments. This identity follows
from the nearest-reference construction and the prescribed
weights. Its management interpretation depends on the
biomass requirements imposed in the model. The agreement
between the analytical expressions and the linear programming
solutions confirms the numerical implementation for the
tested scenarios, while the reported physical quantities
show how each parameter affects the resulting allocation.

\FloatBarrier
\section{Discussion}\label{sec:discussion}

The proposed criterion provides a compromise when the
decision-makers' efficient sets have no common point.
Each set serves as a reference, and the selected decision
minimizes the largest weighted distance to these references.
The formulation also includes the case of compatible
preferences, since its optimal solutions coincide with
the common efficient set whenever that set is nonempty.
The weights and scales determine the meaning of proximity
and therefore form part of the specification of the problem.

The distance criterion also gives a bound on the variation
of each linear objective relative to a nearest efficient
reference. For $x\in X^*$, let $y_p\in E_p$ satisfy
$\|x-y_p\|_\sigma=d_p(x)$. Then, for an objective row $m_{pj}$,
\[
|m_{pj}(x-y_p)|
\le \sum_i |m_{pji}|\,|x_i-y_{pi}|
\le \frac{\tau^*}{\alpha_p}
     \sum_i |m_{pji}|\sigma_i.
\]
Thus, the optimal radius controls objective differences
with respect to the selected references. Estimating losses
relative to ideal objective values would require a separate
comparison with those values.

In the harvest example, the explicit solution clarifies
the respective effects of preferences and physical constraints.
The authority weights determine the selected total harvest,
while the distance scales influence its distribution between
zones. Migration affects the allocation through the biomass
available above the local requirements. For equal weights,
the selected compromise harvests half the largest feasible
total. This relation follows from the geometry of the
reference sets and the chosen selection rule.
Moreover, the opposing objectives $(u,-u)$ make every
feasible allocation efficient for the combined objective
vector. The distance criterion supplies a further basis
for choosing among these allocations.

The aggregate data provide the numerical scale of the example,
while the spatial distribution, migration fractions and
biomass requirements are scenario assumptions. The results
therefore describe the consequences of these assumptions
within a fixed-period model. Their use in a management study
would require local observations and a population model
accounting for recruitment, mortality and changes over time.

For the general problem, the finite linear programming
formulation makes the role of the efficient faces explicit.
Its computational cost depends on the number of faces
and their possible combinations. The allocation example
admits a simpler calculation because its reference sets
reduce to a resource face and the zero-harvest point.
The scalar threshold then yields an optimal selection
and identifies changes caused by active capacity constraints.
The stability estimate complements this formulation by
bounding the effect of approximating the reference sets,
provided that their Hausdorff errors can be controlled.

These results connect the geometry of incompatible efficient
sets with a computable compromise criterion. The finite
formulation gives an exact characterization of the optimal
radius, while the allocation analysis shows how the resulting
decisions depend on weights, scales and resource constraints.
\Needspace{18\baselineskip}
\section{Conclusion}\label{sec:conclusion} We have studied a compromise criterion for decision-makers whose efficient sets may have an empty intersection. The criterion minimizes the largest weighted distance to these sets and recovers their intersection whenever it is nonempty. For linear objectives over a compact polytope, the representation of efficient sets as finite unions of faces gives an exact formulation through a finite family of linear programs. For the harvest-allocation problem, a scalar threshold determines a nearest efficient reference and yields an explicit optimal compromise. The two-zone example shows how the selected allocation depends on the authority weights, distance scales and physical constraints. In particular, migration affects harvest decisions when the resulting biomass availability restricts the feasible allocations. These results provide a direct connection between the geometric formulation and the interpretation of the computed compromise.

\section*{Declaration of AI assistance}

The author used ChatGPT to assist with proofreading, numerical computations, and data provision. The author retains full responsibility for the manuscript, while all key ideas and scientific contributions remain their own.

\bibliographystyle{abbrvnat}
\bibliography{references}

@misc{ices,
  author = {{International Council for the Exploration of the Sea}},
  title = {Stock Assessment Graphs: Northern Norwegian Coastal Cod
           (cod.27.1-2.coastN)},
  year = {2025},
  note = {Assessment key 20955, data for 2024},
  url = {https://standardgraphs.ices.dk/ViewSourceData.aspx?key=20955}
}

@article{Ehrgott2026,
  author  = {Ehrgott, Matthias and K{\"o}ksalan, Murat
             and Kadzi{\'n}ski, Mi{\l}osz and Deb, Kalyanmoy},
  title   = {Fifty years of multi-objective optimization and
             decision-making: From mathematical programming
             to evolutionary computation},
  journal = {European Journal of Operational Research},
  year    = {2026},
  volume  = {330},
  number  = {1},
  pages   = {1--25},
  doi     = {10.1016/j.ejor.2025.06.012}
}

@article{YuZeleny1975,
  author  = {Yu, Po-Lung and Zeleny, Milan},
  title   = {The set of all nondominated solutions in linear cases
             and a multicriteria simplex method},
  journal = {Journal of Mathematical Analysis and Applications},
  year    = {1975},
  volume  = {49},
  number  = {2},
  pages   = {430--468},
  doi     = {10.1016/0022-247X(75)90189-4}
}

@article{Armand1993,
  author  = {Armand, Paul},
  title   = {Finding all maximal efficient faces in
             multiobjective linear programming},
  journal = {Mathematical Programming},
  year    = {1993},
  volume  = {61},
  pages   = {357--375}
}

@article{Sayin1996,
  author  = {Say{\i}n, Serpil},
  title   = {An algorithm based on facial decomposition for finding
             the efficient set in multiple objective linear programming},
  journal = {Operations Research Letters},
  year    = {1996},
  volume  = {19},
  number  = {2},
  pages   = {87--94},
  doi     = {10.1016/0167-6377(95)00046-1}
}

@article{BensonSayin1997,
  author  = {Benson, Harold P. and Say{\i}n, Serpil},
  title   = {Towards finding global representations of the efficient
             set in multiple objective mathematical programming},
  journal = {Naval Research Logistics},
  year    = {1997},
  volume  = {44},
  number  = {1},
  pages   = {47--67},
  doi     = {10.1002/(SICI)1520-6750(199702)44:1<47::AID-NAV3>3.0.CO;2-M}
}

@article{Sayin2000,
  author  = {Say{\i}n, Serpil},
  title   = {Optimizing over the efficient set using
             a top-down search of faces},
  journal = {Operations Research},
  year    = {2000},
  volume  = {48},
  number  = {1},
  pages   = {65--72},
  doi     = {10.1287/opre.48.1.65.12449}
}

@article{TohidiHassasi2018,
  author  = {Tohidi, G. and Hassasi, H.},
  title   = {Adjacency-based local top-down search method for finding
             maximal efficient faces in multiple objective
             linear programming},
  journal = {Naval Research Logistics},
  year    = {2018},
  volume  = {65},
  number  = {3},
  pages   = {203--217},
  doi     = {10.1002/nav.21805}
}

@article{Benson1998,
  author  = {Benson, Harold P.},
  title   = {An outer approximation algorithm for generating all
             efficient extreme points in the outcome set of a
             multiple objective linear programming problem},
  journal = {Journal of Global Optimization},
  year    = {1998},
  volume  = {13},
  number  = {1},
  pages   = {1--24},
  doi     = {10.1023/A:1008215702611}
}

@article{KaciRadjef2022,
  author  = {Kaci, Mustapha and Radjef, Sonia},
  title   = {The set of all the possible compromises of a multi-level
             multi-objective linear programming problem},
  journal = {Croatian Operational Research Review},
  year    = {2022},
  volume  = {13},
  number  = {1},
  pages   = {13--30},
  doi     = {10.17535/crorr.2022.0002}
}

@article{KaciRadjef2023,
  author  = {Kaci, Mustapha and Radjef, Sonia},
  title   = {An adaptive method to solve multilevel multiobjective
             linear programming problems},
  journal = {Operations Research and Decisions},
  year    = {2023},
  volume  = {33},
  number  = {3},
  pages   = {29--44},
  doi     = {10.37190/ord230302}
}

@article{NamNguyenSalinas2012,
  author  = {Nam, Nguyen Mau and Nguyen, Thai An and Salinas, Juan},
  title   = {Applications of convex analysis to the smallest
             intersecting ball problem},
  journal = {Journal of Convex Analysis},
  year    = {2012},
  volume  = {19},
  number  = {2},
  pages   = {497--518},
  url     = {https://journalofconvexanalysis.com/articles/jca19027/}
}

@article{MordukhovichNamVillalobos2013,
  author  = {Mordukhovich, Boris S. and Nam, Nguyen Mau
             and Villalobos, Cristina},
  title   = {The smallest enclosing ball problem and the smallest
             intersecting ball problem: Existence and uniqueness
             of solutions},
  journal = {Optimization Letters},
  year    = {2013},
  volume  = {7},
  number  = {5},
  pages   = {839--853},
  doi     = {10.1007/s11590-012-0483-7}
}

@article{NamVillalobosAn2012,
  author  = {Nam, Nguyen Mau and Villalobos, Maria Cristina
             and An, Nguyen Thai},
  title   = {Minimal time functions and the smallest intersecting
             ball problem with unbounded dynamics},
  journal = {Journal of Optimization Theory and Applications},
  year    = {2012},
  volume  = {154},
  number  = {3},
  pages   = {768--791},
  doi     = {10.1007/s10957-012-0048-z}
}

@article{Patriksson2008,
  author  = {Patriksson, Michael},
  title   = {A survey on the continuous nonlinear resource
             allocation problem},
  journal = {European Journal of Operational Research},
  year    = {2008},
  volume  = {185},
  number  = {1},
  pages   = {1--46},
  doi     = {10.1016/j.ejor.2006.12.006}
}

@article{Elbetch2025,
  author  = {Elbetch, Bilel and Moussaoui, Ali and Auger, Pierre},
  title   = {Increase maximum economic yield in a patchy environment},
  journal = {Journal of Mathematical Biology},
  year    = {2025},
  volume  = {90},
  number  = {1},
  pages   = {14},
  doi     = {10.1007/s00285-024-02178-6},
  note    = {Article 14}
}

@article{Kaci2025,
  author  = {Kaci, Mustapha},
  title   = {Optimizing multilevel decentralized problems
             with a novel application of Gabasov's adaptive
             method: highlighting applications in COVID-19
             vaccine distribution},
  journal = {RAIRO - Operations Research},
  volume  = {59},
  pages   = {803--833},
  year    = {2025},
  doi     = {10.1051/ro/2025010}
}
\end{document}